\documentclass[11pt]{amsart}

\usepackage{microtype}
\usepackage{amsmath, amssymb, amsthm}
\usepackage[hidelinks]{hyperref}
\usepackage[parfill]{parskip}        
\usepackage[margin=1.16in]{geometry} 
\usepackage{tikz}
\usepackage{tikz-cd}
\usetikzlibrary{matrix}

\usepackage{mathtools}

\title{Aperiodicity, rotational tiling spaces and topological space groups}\keywords{} \subjclass{}

\keywords{Aperiodic tilings, tiling spaces, rotations, space groups, crystallographic groups, shape theory, group cohomology.} 

\date{\today}

\author{John Hunton}
\address{Department of Mathematical Sciences, Mathematical Sciences \& Computer Science Building, Durham University, Upper Mountjoy Campus, Stockton Road, Durham, DH1 3LE, UK}
\email{john.hunton@durham.ac.uk}
\urladdr{http://www.dur.ac.uk/john.hunton}

\author[James Walton]{James J.\ Walton}
\thanks{Research of the second author supported by EPSRC grant EP/R013691/1.}
\address{School of Mathematical Sciences, University of Nottingham, University Park, Nottingham, NG7 2RD, UK}
\email{James.Walton@nottingham.ac.uk}
\urladdr{https://www.nottingham.ac.uk/mathematics/people/james.walton}

\newcommand{\R}{\mathbb{R}}
\newcommand{\Q}{\mathbb{Q}}
\newcommand{\Z}{\mathbb{Z}}
\newcommand{\N}{\mathbb{N}}
\newcommand{\EE}{\mathcal{E}}
\newcommand{\Om}{\Omega}
\newcommand{\T}{{\mathbb T}}
\newcommand{\F}{{\mathbb F}_2}
\newcommand{\FF}{{\mathcal F}}
\newcommand{\SO}{{\mathrm{SO}}}
\newcommand{\Sp}{{\mathrm{Spin}}}
\newcommand{\Or}{{\Omega_r}}
\newcommand{\Ot}{{\Omega_t}}
\newcommand{\GG}{{\widetilde G}}
\newcommand{\ph}{\pi^{\mathrm{pro}}} 
\newcommand{\sg}{\Gamma^\mathrm{pro}_{\phantom{+}}} 
\newcommand{\psg}{\Gamma^\mathrm{pro}_+} 
\newcommand{\lra}{\longrightarrow}
\newcommand{\lla}{\longleftarrow}
\newcommand{\il}{\varprojlim}
\newcommand{\asg}{topological space group} 
\newcommand{\spg}{space pro-group} 

\newcommand{\MM}{\mathbb{M}}
\newcommand{\XX}{\mathbb{X}}
\newcommand{\PP}{\mathbb{P}}

\newtheorem{theorem}{Theorem}[section]
\newtheorem{proposition}[theorem]{Proposition}
\newtheorem{lemma}[theorem]{Lemma}
\newtheorem{corollary}[theorem]{Corollary}

\theoremstyle{definition}
\newtheorem{definition}[theorem]{Definition}
\newtheorem{example}[theorem]{Example}
\newtheorem{examples}[theorem]{Examples}
\newtheorem{remark}[theorem]{Remark}

\numberwithin{equation}{section} 
\numberwithin{figure}{section}   

\begin{document}

\begin{abstract}
We study the rotational structures of aperiodic tilings in Euclidean space of arbitrary dimension using topological methods. Classical topological approaches to the study of aperiodic patterns have largely concentrated just on translational structures, studying an associated space, the continuous hull, here denoted $\Om_t$. In this article we consider two further spaces $\Om_r$ and $\Om_G$ (the rotational hulls) which capture the full rigid motion properties of the underlying patterns. The rotational hull $\Om_r$ is shown to be a matchbox manifold which contains $\Om_t$ as a sub-matchbox manifold. We develop new S-MLD invariants derived from the homotopical and cohomological properties of these spaces demonstrating their computational as well as theoretical utility. We compute these invariants for a variety of examples, including a class of 3-dimensional aperiodic patterns, as well as for the space of periodic tessellations of $\R^3$ by unit cubes. We show that the classical space group of symmetries of a periodic pattern may be recovered as the fundamental group of our space $\Om_G$. Similarly, for those patterns associated to quasicrystals, the crystallographers' aperiodic space group may be recovered as a quotient of our fundamental invariant. 
\end{abstract}

\maketitle

\section{Introduction}
Space groups, also known as crystallographic groups or Bieberbach groups, capture the symmetries of periodic patterns or tilings in Euclidean space $\mathbb{R}^d$. The space group of a particular pattern $T\subset\mathbb{R}^d$ is the subgroup of the full isometry group of $\mathbb{R}^d$ that fixes $T$; such groups have been classical objects both of study and application in Mathematics, Physics and Chemistry since at least the $19^{\rm th}$ century.

In this article we consider {\em aperiodically ordered\/} patterns, an infinite class of highly structured but non-periodic patterns in $\mathbb{R}^d$. This class includes, as special cases, various well-known examples such as the Penrose Tilings, as well as the objects used to model quasicrystals, materials only discovered at the end of the $20^{\rm th}$ century \cite{ShBlGrCa84}. For those patterns used to model quasicrystals there is a well developed analogue to the space group of a periodic pattern \cite{dWJaJa81, Mer92, Jan88, RoWrMe88, RabFis03} which captures information about their rotational and translational structures relating to their diffraction images. In this article we use a topological approach to define and compute new algebraic invariants for aperiodic patterns. Our first invariant, in the case of the patterns modelling quasicrystals, has the crystallographers' aperiodic space group as, typically, a strict quotient.

The association of a `space group' to a non-periodic pattern may seem at first sight an oxymoron: the traditional association of a group to a periodic pattern is an algebraic encoding of the precise symmetries it enjoys, but it is the very nature of an aperiodic pattern that it is notably short of symmetries, at least considered translationally. However, methods from Topology and Dynamics (see, for example \cite{BaGr}, or the collection of surveys \cite{AObook}) have long proved effective in capturing various structural properties in the absence of exact symmetries. Topology in particular has  come into the study of aperiodic patterns as follows. To any periodic or aperiodic tiling $T$, a topological space $\Om$ known as the {\em tiling space\/} or {\em continuous hull}, is constructed from the set of translated images of $T$. The properties of $\Om$, particularly as seen by tools such as \v{C}ech cohomology or $K$-theory, reveal many key aspects of the translational structure of $T$: see \cite{HunMFO17, SadBook08} for a discussion and brief survey of some of the more significant results of this  approach. For a periodic pattern, $\Om$ is just a $d$-torus, and for pertinent reasons this should be thought of as the classifying space of the group $\Z^d$ of translational symmetries of $T$; for an aperiodic example, $\Om$ may be seen as a natural generalisation in that it can be realised as the classifying space of a certain associated translation {\em groupoid}, an aspect to which we shall return in a further article.

Nevertheless, with a few notable exceptions, such as \cite{BSJ91, Mal15, Rad94, Ran07, Sta15, Wal17rot}, and see also \cite{BaGr}, the topological study to date has largely been confined to the analogue of the translational symmetries for aperiodic patterns. This is perhaps surprising as many of the most interesting examples display apparent strong rotational or reflective organisation. Indeed, although it is the translational structure that determines in what way the pattern is diffractive, it is precisely the rotational structure that gives rise to the rotational properties of any associated diffraction pattern, the properties that first alerted researchers to consider them as models for quasicrystals.

Considering only the translational structure misses a good deal. For example, in the case of periodic tilings in the plane, it is well known that there are 17 \lq wallpaper groups\rq\ that can act -- but the subgroup of translations in all cases is just free abelian of rank 2. In three dimensions, there are over 200\footnote{219, or 230 if chiral pairs are distinguished.} \lq crystallographic groups\rq, but again all translation subgroups of these are isomorphic. In general $\mathbb{R}^d$, the rigid symmetries of a periodic pattern are captured by its space group $\Gamma$ which can always be described as an extension 
\[
0\lra\Z^d\lra \Gamma\lra G\to 1
\]
of the translation subgroup by the finite {\em point group\/} $G$. For given $\Z^d$ and $G$ there is usually more than one such possible extension. For aperiodic patterns modelling quasicrystals, the corresponding space group is again finite by free abelian, though the rank of the free abelian subgroup is larger than the dimension $d$ of the quasicrystal. 

In general, sufficiently regular (for example, repetitive) aperiodic patterns have an analogous \lq point group\rq\, which likewise captures richer structure about the pattern than is capable with the purely translationally defined tiling space alone; in the aperiodic case this group can even be infinite, for example it is the group $\mathrm{O}(2)$ in the case of the Pinwheel tiling \cite{Rad94}. The notion of point group is already documented in the literature, see for example \cite{BaaGri12}, but we discuss it in detail for general patterns in Section \ref{sec: rotations}, defining precisely the class of aperiodic tilings considered in Section \ref{sec: hull}. This includes all repetitive tilings (Definition \ref{def: repetitive}), which are shown to always have well-defined point groups in Proposition \ref{prop: repetitive => pg}.

In this article we develop the topological approach beyond  translational issues, we study this richer structure for aperiodic patterns  also by way of an associated topological space, in fact by two such spaces. The first of these is the {\em rotational hull\/} $\Om_r$. Just as $\Om$, the continuous, or what we shall now refer to as the {\em translational}, hull, denoting it by $\Om_t$, can be considered as a certain completion of the space of translates of the pattern, the rotational hull $\Om_r$ is the corresponding completion of the set of all Euclidean motions of $T$. The  space $\Om_r$ has been considered before, but by and large only for 2-dimensional patterns \cite{BDHS10, SadBook08, Wal17rot}. Like $\Om_t$, the space $\Om_r$ is also a matchbox manifold (Proposition \ref{OmrMM}), but now of dimension $d+d(d-1)/2$, where $d$ is the dimension of the Euclidean space in which the pattern lives. ($\Om_t$ is a matchbox manifold of dimension $d$.)

Our second space, which we denote $\Om_G$, is the Borel construction (or `homotopy quotient' of $\Om_t$ by $G$) arising from the action of the aperiodic point group $G$ on $\Om_t$; this space lies at the heart of our topological analogue of the space group. We formally define the spaces $\Om_t$ and $\Om_r$ in Section \ref{sec: hull} and $\Om_G$ in Section \ref{sec: useful}, where we also relate these three spaces.

The main results of this paper provide an analysis of the algebraic topology of $\Om_r$ and $\Om_G$, and this allows us to define algebraic objects associated to an aperiodic pattern that are invariant under a standard, natural notion of equivalence (specifically, under S-MLD equivalence, see \cite{BSJ91}). We consider both homotopical and cohomological viewpoints and  present a framework for computation for patterns in any dimension of space $\R^d$, $d\geqslant2$, Section \ref{sec: useful}. The spaces $\Om_r$ and $\Om_G$ are of course closely related, and hence so are their topological invariants. The homotopy theory of $\Om_G$ appears to be the more fundamental, but the space $\Om_r$ has the practical advantage of being cohomologically finite dimensional, and so may be seen as a useful staging post on the way to mining the richer information in $\Om_G$.

The computation of any of these invariants  even for the translational hulls is frequently difficult, especially as the dimension $d$ increases: computation for the rotational analogues is typically more complex still. Nevertheless, we are able to present a number of worked examples to demonstrate the practicality of the machinery we develop. 

In fact, there is some merit even in applying our approach to periodic examples, and in Section \ref{sec: cube} we give a complete computation for the cohomology   $H^*(\Om_r;\Z)$ in the case of the periodic  tessellation of $\R^3$ by unit cubes. This is effectively a computation of the cohomology of the 6-manifold of configurations of the  cubical lattice in $\R^3$, and may be of independent interest. This computation then provides the foundation for our final result, the computation of the integral cohomology of the rotational hull of a class of 3-dimensional aperiodic tilings based on decorated cubes, Section \ref{sec: Fib3}.

More can be said in the case of rational cohomology, which gives a less subtle but more easily computable invariant. In Section \ref{sec: rational} we give for tilings $T$ in all dimensions $d\geqslant 2$ a complete description of the rational cohomology $H^*(\Om_r;\Q)$ in terms of the aperiodic point group $G$ and its action on $H^*(\Om_t;\Q)$. These calculations determine the ranks of $H^*(\Om_r;\Q)$ (or equivalently $H^*(\Om_r,\R)$) in terms of invariants from the translational setting, for which there are well-established tools of computation. We note that the cohomology groups $H^*(\Om_t;F)$, for $F = \Q$ or $\R$, have been used in the past for trace purposes, most notably in Bellissard's Gap Labelling Theorem \cite{BBG06}. See also \cite[Section 9]{AndPut98}, and \cite{Sad11}, where the top degree $H^d(\Om_t,\Q)$ plays a distinguished role. It seems natural to exploit the extra structure provided by the action of rotational symmetry in such constructions. In Section \ref{sec: top degree} we identify the top degree rational cohomology of $\Om_r$ with the subgroup of elements of the top degree rational cohomology of $\Om_t$ invariant under the action of the point group $G$, from which there exists a natural trace map induced by patch frequencies.

The integral cohomology requires deeper input. The case of integral cohomology of planar tilings has already been covered by the second author \cite{Wal17rot}, and in Section \ref{sec: planar} we recover the final descriptions of that work via the machinery of the current article. In higher dimensions, full computations are difficult. Several of the results in Section \ref{sec: cohomology} provide the first steps to effective calculation of these invariants, which we hope will provide the foundations for more powerful tools in the future. One technique which has proved to be useful is the replacement of the space $\Om_r$ (and related spaces, such as $\Om_t$) with shape equivalent \cite{MarSeg82} but less pathological cellular approximations, which have isomorphic cohomology.

The homotopy groups associated to aperiodic tilings have been less well studied, though we note the pioneering work of Geller and Propp \cite{GelPro95} in the translational context. As already indicated, although $\Om_r$ and $\Om_G$ are closely related, here it would seem the latter which is the more fundamental.

For a periodic tiling $T$, the fundamental group $\pi_1(\Omega_G)$ is naturally isomorphic to the  group $\Gamma$ of symmetries of $T$, i.e., its space group (Corollary \ref{cor: per pi1 Om_G}). There is a canonical cover of $\Gamma$ (for $d=2$ it is a $\Z$ cover which keeps track of winding number information of symmetries, and for $d > 2$ it is a $\Z/2$ cover) which is realised (Corollary \ref{cor: per pi1 Om_r}) by the fundamental group of  $\Omega_r$. For periodic tilings, all information required to reconstruct a tiling (up to a natural notion of locally defined redecoration) is contained in its space group. Thus the homotopical study of $\Om_G$, and less directly $\Omega_r$, is tightly linked with the classical study of such tilings, giving in turn a further context in which to understand what the cohomology of $\Omega_r$ tells us about the tiling. We note in particular the article of Hiller \cite{Hil86} which uses related cohomological tools to classify the periodic space groups.

Although there is no direct analogue of the space group for a general aperiodic tiling, the topological spaces $\Omega_G$ and $\Omega_r$ are still defined. Because of the pathological nature of these spaces, the classical homotopy groups are not well suited to them, but Shape Theory \cite{MarSeg01} provides appropriate replacements via the \emph{shape homotopy groups} \cite{ClHu12, MarSeg82}. We thus introduce our fundamental invariant, our  `\asg' of a pattern as the shape fundamental group of $\Om_G$. This is a natural extension of the notion of space group for a periodic pattern, and in the case of tilings modelling quasicrystals, the crystallographers' aperiodic space group can be derived from the \asg, but the latter appears to be a richer invariant (and one retains more information still by using the `topological space pro-group', Definition \ref{def: top sp gp}). As in the periodic case the shape fundamental pro-group of $\Om_r$ corresponds to an associated cover of the \asg. 

We note that the \asg \ is an invariant of the original tiling $T$, rather than the space $\Om_G$. Indeed, a priori, there may be a base point dependence for aperiodic $T$, although we do not currently know of specific examples for which this is the case.

A full description of the topological space group would seem infeasible in all but the simplest examples (such as those given by products of one-dimensional tilings). However, invariants of the topological space pro-group could be more accessible yet still contain rich information. For example, one may apply the functor $\hom(-,G)$ for a finite (and non-abelian) group $G$ and take the corresponding direct limit. Such invariants were considered by Sadun in \cite[Section 4]{Sad:PEC} and have been calculated by G\"{a}hler for some one- and two-dimensional substitution tilings in the translational setting. Even with relatively small groups $G$, these invariants are often capable of distinguishing examples with isomorphic cohomology.

The \asg\ thus defined is a new invariant associated to tilings, which is invariant under S-MLD equivalence of tilings. Such invariance, together with results relating the shape homotopy groups of $\Om_t$, $\Om_r$, $\Om_G$, are presented in Section \ref{sec: homotopy}, where we also relate the \asg\ to the classical space groups of periodic tilings and quasicrystals. Section  \ref{sec: homotopy} also contains  computations and descriptions for several families of aperiodic examples.

\subsection*{Acknowledgements}
The authors thank Michael Baake and Franz G\"ahler for helpful discussions concerning aperiodic space groups for quasicrystals. We also thank the anonymous referee for their valuable suggestions.


\section{Patterns, tilings and their hulls} \label{sec: hull}
For the purposes of this article, and for simplicity, the patterns we consider are tilings of $\R^d$, $d\geqslant2$, by compact $d$-dimensional polyhedral subsets. The restriction to polyhedral tilings is not an important one. Other patterns -- represented for example by labelled point sets, or tilings of fractal or possibly overlapping tiles -- can, given reasonable restrictions, always be represented by polyhedral tilings that are equivalent. More precisely, such a pattern can always be represented by a polyhedral tiling which is S-MLD equivalent to it (see \cite{BSJ91} and Definition \ref{SMLDdef} below for details of this equivalence relation).

We follow the standard set-up for discussing aperiodic tilings in $\R^d$, and briefly introduce the necessary concepts for the new reader here. Further details may be found, for example, in \cite{BaGr} or \cite{SadBook08} where there is extended commentary on the underlying ideas.

\begin{definition}
A tiling $T$ in $\R^d$ is a cover of $\R^d$ by compact $d$-dimensional polyhedral subsets, called \emph{tiles}, meeting full face to full face, and only ever on boundaries. To distinguish  tiles of the same geometric shape further, it is sometimes convenient to also allow each to carry a `label' or `colour'. We assume that each tile is congruent (with matching labels), by translation, to one of a finite set of (labelled) polyhedra, the {\em prototiles}.
\end{definition}

For example, the `infinite chess board' has two prototiles: a black unit square and a white unit square. The Penrose tiling has 10 prototiles, each congruent (by rotation)  to a unit sided rhombus with either a $\pi/5$ angled corner, or a $2\pi/5$ angled one, each rhombus occurring in one of 5 possible rotations.

\begin{definition}
A \emph{patch} of a tiling $T$ is a finite selection $P \subset T$ of tiles from $T$. The patch of tiles intersecting a closed Euclidean ball of radius $r$ at $x \in \R^d$ is called the \emph{$r$-patch} centred at $x$. 
\end{definition}

With these definitions, the objects we study include most of the standard examples of aperiodic tilings studied in the literature. Under the conditions stated, the tilings necessarily have {\em translational finite local complexity}, that is, for any given radius $r$, there are up to translation only a finite number of $r$-patches; we call such a tiling an \emph{FLC tiling}. Thus we exclude, for example, the Pinwheel tiling \cite{Rad94}: there is good reason for this in that the Pinwheel, whose rotational structure has been well studied already \cite{BDHS10, FrWhWh14}, has an infinite point group; we are concerned here with tilings with finite rotational structure, as defined in the next section.

We introduce the first two tiling spaces considered in this paper. Both use a metric on sets of tilings; in fact there is a considerable choice in the actual metric used and it is the topology they define that really matters, but the key underlying concept is the following. Loosely speaking, two tilings are considered as close if, up to a small perturbation of either, the two agree to a large radius about the origin in $\R^d$. For the definitions below, and under our current assumptions about our tilings, it will be enough to take as a `perturbation' a rigid motion. Thus we take two tilings as close if they agree to a large radius about the origin after a small translation followed by a small rotation of either of them. One can easily make this geometric structure precise by introducing a `tiling metric' \cite{SadBook08} on a given suitable collection of tilings, or a uniformity \cite{Wal17pe}, which is more canonical and more easily verified as providing the required geometric structure.

\begin{definition}
Given a tiling $T$,  say that another $T'$ is \emph{locally indistinguishable} from it if every patch of tiles of $T'$ appears, up to translation, in $T$. We define the \emph{translational hull} or {\em tiling space\/} $\Om_t$ of $T$ to be the topological space of tilings locally indistinguishable from $T$, taken with the topology described above. 
\end{definition}

In fact, the tiling metric or uniformity also provides a notion of a Cauchy sequence of tilings. One may show that $\Omega_t \cong \overline{T + \R^d}$; that is, the translational hull is homeomorphic to the completion of the space of translates of the tiling $T$.

The space $\Om_t$ encodes information related to the translational structure of $T$ topologically. It has been widely studied in the literature, see \cite{SadBook08} for an introduction. In this paper we wish to consider a closely related space which also captures rotational aspects of the tiling. This is done by a natural modification of the above definition: local indistinguishability only allows for comparison of patches by translations. By allowing general rigid motions, we define the rotational hull:

\begin{definition}
The \emph{rotational hull} or {\em rotational tiling space\/} $\Om_r$ of $T$ is defined as the space of tilings whose finite patches all appear, up to rigid motion, in $T$.
\end{definition}

Here, and throughout, `rigid motion' means  an orientation preserving isometry of $\R^d$, for reasons elaborated on below in Remark \ref{rem: reversing}.

Finite local complexity allows us to alternatively define $\Om_r$ as the space of rotates of elements from $\Om_t$:

\begin{proposition} \label{prop: FLC => global rotate}
We have that $T_1 \in \Om_r$ if and only if $T_1 = g(T_2)$, for some $g \in \SO(d)$ and $T_2 \in \Om_t$, that is, $\Om_r = \SO(d) \cdot \Om_t$.
\end{proposition}

\begin{proof}
If $T_2 \in \Om_t$ then every finite patch in $T_2$ appears in $T$ up to translation, and hence every finite patch of $T_1 = g(T_2)$ appears in $T$ up to rigid motion, so $\SO(d)\cdot \Om_t \subseteq \Om_r$. Suppose then that $T_1 \in \Om_r$. Let $P_n$ denote the patch of radius $n$ centred at the origin in $T_1$. Since these patches appear in $T$ up to rigid motion, there exist $g_n \in \SO(d)$ for which $g_n(P_n)$ appears in $T$ up to translation, for all $n \in \N$. By FLC, there are only finitely many such $g_n$, so $g_n = h$ for infinitely many $n$, for some $h \in \SO(d)$. Given $g_n = h$, we have that $h(P_j) \subseteq h(P_n)$ also appears in $T$ up to translation for any $j \leq n$, so we may take $g_n = h$ for all $n \in \N$. Hence, every $n$-patch centred at the origin in $h(T_1)$ appears in $T$, up to translation. Since every finite patch of $h(T_1)$ is eventually contained in such a patch, we see that $h(T_1) \in \Om_t$. So $T_1 = g(T_2)$, where we take $T_2 = h(T_1) \in \Om_t$ and $g = h^{-1} \in \SO(d)$, as required.
\end{proof}

As for the translational hull, the rotational hull can also be described as a completion (for more details, see \cite{BDHS10}): $\Om_r$ is the completion of the space of rigid motions of $T$. This allows for a shorter proof of the above, noting that $\SO(d) \cdot \Om_t$ is a compact subset of $\Om_r$ containing the Euclidean orbit of $T$, and so is all of $\Om_r$.

\begin{remark}
It follows immediately from these definitions that $\Om_t$ has a natural $\R^d$ action, by translation, and $\Om_r$ has a natural action by the (positive) Euclidean group of rigid motions: if $\Phi$ is a translation, respectively a rigid motion, in $\R^d$ and $T'\in\Om_t$, respectively $\Om_r$, then $\Phi(T')$ is also an element of the corresponding tiling space. In particular, $\Om_r$ has an action by $\SO(d)$. It may readily be checked that these actions are continuous from the definition of the underlying topology.
\end{remark}

The issue of when two tilings should be considered `equivalent' is an important one. A key concept is that of {\em Mutually Locally Derived\/} (MLD) equivalence, and, when we consider rotational structures, its analogue {\em S-MLD\/} equivalence \cite{BSJ91}. 

\begin{definition}
Given tilings $T_1$ and $T_2$, we say that $T_2$ is \emph{locally derivable} from $T_1$ if there exists some $r>0$ for which, whenever $\Phi$ is a translation which identifies the $r$-patch at $x$ of $T_1$ with that at $\Phi(x)$, then $\Phi$ also identifies the $1$-patch at $x$ of $T_2$ to that at $\Phi(x)$.
\end{definition}

\begin{definition}\label{SMLDdef}
Given tilings $T_1$ and $T_2$, we say that $T_2$ is \emph{S-locally derivable} from $T_1$ if there exists some $r>0$ for which, whenever $\Phi$ is a rigid motion which identifies the $r$-patch at $x$ of $T_1$ with that at $\Phi(x)$, then $\Phi$ also identifies the $1$-patch at $x$ of $T_2$ to that at $\Phi(x)$.
\end{definition}

The choice of radius $1$ in these definitions is arbitrary: the point is that to decide how the pattern of $T_2$ is tiled locally about a point $x \in \R^d$, one only needs to know the decoration of $T_1$ to radius $r$ centred at $x$, up to rigid motion. Thus the local derivation is encoded by a rule for redecorating $T_1$ to get $T_2$ which is locally defined, and respecting rotational symmetries for an $S$-local derivation. 

\begin{definition}
If $T_2$ is locally derivable, respectively S-locally derivable, from $T_1$ and vice versa, then we call $T_1$ and $T_2$ \emph{MLD}, respectively \emph{S-MLD}. These are equivalence relations on sets of tilings. 
\end{definition}

These equivalence relations can be generalised to patterns which are not necessarily polyhedral tilings, and one could incorporate orientation reversing symmetries too if desired.

\begin{remark}
If $T_2$ is locally derivable from $T_1$ then there is an induced map from the translational hull for $T_1$ to that for $T_2$. It follows that if $T_1$ and $T_2$ are MLD equivalent, then their translational hulls are homeomorphic, in fact with homeomorphism commuting with the translation action by $\R^d$. Similarly, if $T_1$ and $T_2$ are S-MLD equivalent, then their rotational hulls are homeomorphic with homeomorphism commuting with the action of rigid motions, in particular by rotations in $\SO(d)$.
\end{remark}

Most tilings currently of interest in the field of Aperiodic Order are \emph{repetitive}.

\begin{definition} \label{def: repetitive}
A tiling $T$ is \emph{repetitive} if for each $r>0$ there is an $R_r>0$ such that every $r$-patch can be found, up to translation, within distance $R_r $ of every point $x \in \R^d$.
\end{definition}

For a repetitive tiling $T$, and $T'$ locally indistinguishable from $T$, it is easily shown that $T'$ has the same set of finite patches, up to translation (that is, $T$ is also locally indistinguishable from $T'$). Hence, repetitivity ensures that if $\Om$ is the hull of $T$ (translational or rotational), $T' \in \Om$ and $\Om'$ is the hull of $T'$, then $\Om = \Om'$. This fails if $T$ is not repetitive. Repetitivity is equivalent to minimality of the dynamical system $(\Om_t,\R^d)$ \cite{AObook}.

\begin{remark}\label{rem: reversing}
The symmetry structure for tilings considered throughout this paper is mostly confined to orientation preserving symmetries; the reader may reasonably wonder why we do not consider orientation reversing symmetries of our tilings. Certainly the existence or not of orientation reversing symmetries of a given pattern is important, and extending some of the constructions here to orientation reversing symmetries is essentially straightforward, but provides limited additional insight. One of our central focuses here is an analysis of the topology of the rotational hull $\Om_r$. Suppose we define $\Om_s$ as the corresponding completion but now of all Euclidean motions, orientation reversing as well as preserving. For a repetitive tiling, if every finite patch which occurs in $T$ also has its mirror image occurring in $T$ (up to some rigid motion) then $\Om_r = \Om_s$, which additionally carries a $\Z/2$ action corresponding to a reflection. If finite patches do not have mirror images in the tiling  then $\Om_s$ is homeomorphic to the disjoint union of two copies of $\Om_r$ (and the additional $\Z/2$ action merely swaps the two components), so in either case it suffices to study this space $\Om_r$ alone.
\end{remark}

\subsection{Global structure of the rotational hull} \label{sec: matchbox}
In the translational setting, it is well known that $\Om_t$ is an orientable \emph{matchbox manifold}, see \cite{AaHaOv91, CandelConlon,  MooreSchochet}. In brief, a matchbox manifold is a continuum foliated with Euclidean leaves with totally disconnected local transversals. If we can take for each transversal just a single point, this is an ordinary manifold, as is the case for periodic tilings, but aperiodicity forces us to need the richer structure; the leaves of the foliation for $\Om_t$ are the path components, given by the $d$-dimensional translational orbits. A similar result, Proposition \ref{OmrMM}, holds for $\Om_r$: the rotational hull is a matchbox manifold of dimension $d+d(d-1)/2$ whose leaves are the orbits of tilings under \emph{rigid motion}. We will explain in this section how charts are constructed to give this local product structure, but first we begin by recalling the definition of a matchbox manifold. We follow \cite{CandelConlon, ClHu13, ClHuLu14, ClHuLu19, MooreSchochet} where the reader may find further discussion of this concept.

\begin{definition}\label{MatMan}
A matchbox manifold of dimension $n$ is a continuum $\MM$ satisfying the following conditions. 
\begin{itemize}
\item {\em Charts.} There is a compact, separable, totally disconnected metric space $\XX$, and for each $x\in\MM$ a compact subspace $X_x\subset\XX$, an open set $U_x\subset\MM$, and a homeomorphism $\phi_x\colon \overline{U}_x\to [-1,1]^n\times X_x$ on its closure $\overline{U}_x$ in $\MM$ such that $\phi_x(x)=(0,w_x)$ for some $w_x\in \mbox{int}(X_x)$. Moreover, we assume $\phi_x$ may be extended to a homeomorphism $\hat\phi_x\colon \hat{U}_x\to [-2,2]^n\times X_x$ for some open set $\hat{U}_x\subset\MM$ with $\overline{U}_x\subset\hat{U}_x$.
\item {\em Plaques.}  Let $\pi_x\colon\overline{U}_x\to X_x$ denote the composite of $\phi_x$ followed by projection onto the second factor, and for $z\in \overline{U}_x$ define the plaque $\PP_x(z)$ through $z$ for the chart $\phi_x$ as $\pi_x^{-1}(\pi_x(z))\subset \overline{U}_x$. Each plaque $\PP_x(z)$ is given the topology such that the restriction $\phi_x\colon\PP_x(z)\to [-1,1]^n\times\{\pi_x(z)\}$ is a homeomorphism.
\item{\em Compatibilities.} For two charts $(U_x,\phi_x)$ and $(U_y,\phi_y)$, any intersection $\PP_x(z)\cap \PP_y(z')$ is open in each plaque,  each plaque $\PP_x(z)$ in $U_x$ meets at most one plaque in $\overline{U}_y$, and for each $z\in \PP_x(z)\cap \PP_y(z)$  there is a connected  set $W\subset \PP_x(z)\cap \PP_y(z)$, open in the plaque topology, containing $z$, and such that 
$$\phi_y\circ\phi_x^{-1}\colon \PP_x(z)\cap W\to \PP_y(z)\cap W$$
is a smooth diffeomorphism from 
$\phi_x(\PP_x(z)\cap W)$ to $\phi_y(\PP_y(z)\cap W)$.
\end{itemize}
\end{definition} 

Our result, Proposition \ref{OmrMM}, holds for tilings that only have \emph{Euclidean finite local complexity}, that is, tilings for which, for each $r > 0$, there are only finitely many $r$-patches up to rigid motion (rather than the more restrictive condition of agreeing up to translation). Therefore, in this section we temporarily weaken our usual FLC condition to Euclidean FLC.

We choose a `puncture' (an interior point) for each tile, so that isometric tiles are punctured identically. We assume this is done so that all possible symmetries of a tile preserve the puncture. This may be achieved, for example, by taking tiles to be convex polytopes (which may always be done, up to S-MLD equivalence) and taking punctures to be centres of mass of tiles. By Euclidean FLC, there are only finitely many tiles in $T$ up to rigid motion. We may thus choose a finite set of representative \emph{prototiles} $\mathcal{P}$ satisfying the following:
\begin{itemize}
	\item each $t \in \mathcal{P}$ is a tile appearing in $T$, up to rigid motion;
	\item every tile of $T$ is a rigid motion of some $t \in \mathcal{P}$;
	\item each $t \in \mathcal{P}$ has puncture over the origin.
\end{itemize}
Of course, we could assert that each tile of $T$ is the rigid motion of a \emph{unique} tile from $\mathcal{P}$, but for analogy later with the translational setting we allow multiple rotates of tiles in $\mathcal{P}$, which will not affect our arguments.

We define the (\emph{canonical}) \emph{transversal} as
\[
\Xi \coloneqq \{T' \in \Om_r \mid t \text{ appears in } T\text{, for } t \in \mathcal{P}\}.
\]
That is, we consider the collection $\Xi \subset \Om_r$ of tilings with a puncture at the origin, so that the tile containing the origin is oriented identically to a prototile in $\mathcal{P}$.

\begin{lemma}
The transversal $\Xi$ is a compact, separable, totally disconnected space.
\end{lemma}

\begin{proof}
The elements of $\Xi$ may be identified with sequences
\[
P_0 \subset P_1 \subset P_2 \subset P_3 \subset P_4 \subset \cdots,
\]
where $P_0 = \{t\}$ for some $t \in \mathcal{P}$, and for higher $n \in \N$ each $P_n$ is an $n$-patch, centred at the origin, which appears in $T$ up to rigid motion and extends $P_{n-1}$. Since the $P_n$ cover all of $\R^d$, such a sequence defines and is defined by a tiling of $\Xi$.

The prototile set $\mathcal{P}$ is finite, and by Euclidean FLC there are only finitely many ways of extending a finite patch to another. Hence there are only finitely many possibilities for each $P_i$. Two such sequences are close in the tiling topology if and only if they agree for large $n$. So, as a subspace of $\Om_r$ with the tiling topology, $\Xi$ is an inverse limit of discrete, finite spaces.
\end{proof}

Small but distinct rigid motions move tilings of $\Xi$ to distinct tilings in $\Om_r$. To see this, first let $\epsilon > 0$ be such that all punctures are distance greater than $\epsilon$ from the boundaries of tiles. Choose a small neighbourhood $U$ of the origin in $\SO(d)$, homeomorphic to the closed ball in $\R^{d(d-1)/2}$, so that for all $u_1, u_2 \in U$ and $t_1, t_2 \in \mathcal{P}$, if $u_1(t_1) = u_2(t_2)$ then $u_1 = u_2$. Equivalently, for any $t \in \mathcal{P}$ and non-trivial $g \in U \circ U^{-1}$, we have that $g(t) \notin \mathcal{P}$. Since $\SO(d)$ is a Lie group of dimension $d(d-1)/2$, and sufficiently small non-trivial rotates of prototiles do not appear in $\mathcal{P}$ (since $\mathcal{P}$ is a finite set of convex polytopes), such a $U \subset \SO(d)$ exists.

Define $B = B_\epsilon \times U$, where $B_\epsilon$ is the set of translations of $\R^d$ of norm at most $\epsilon$, which we may identify with the unit ball of $\R^d$. So we may identify $B$ with a Euclidean ball of dimension $d+d(d-1)/2$. We associate each $(\tau,g) \in B$ with the rigid motion $\tau \circ g$ and define
\[
f \colon B \times \Xi \to \Om_r, \ f(b,Y) \coloneqq b(Y)
\]
Let $X$ denote the image $f(B \times \Xi)$.

\begin{lemma}
The map $f \colon B \times \Xi \to X$ is a homeomorphism.
\end{lemma}

\begin{proof}
Continuity of $f$ follows from the definition of the tiling topology. By definition $f$ is surjective onto $X$. To see that $f$ is injective, suppose that $g_1(T_1) + x_1 = g_2(T_2) + x_2$ for some $x_1$, $x_2 \in \R^d$ with norm at most $\epsilon$, rotations $g_i \in U$ and $T_i \in \Xi$. Each $g_i(T_i)$ is a tiling with puncture over the origin. Since $x_1$ and $x_2$ have norm at most $\epsilon$, translates of the same tiles contain the origin after translation by $x_1$ and $x_2$, which thus now have punctures at $x_1$ and $x_2$. It follows that $x_1 = x_2$, since we have assumed that punctures can be uniquely determined from the geometry of the tiles. By translating back by $-x_1=-x_2$, we see that $g_1(T_1) = g_2(T_2)$. Since each $T_i \in \Xi$ contains a tile of $\mathcal{P}$ at the origin, by the definition of $U$ we have that $g_1 = g_2$. Similarly, by applying $g_1^{-1}$, it follows that $T_1 = T_2$ and $f$ is injective, as required. Since $B \times \Xi$ is compact and $\Om_r$ is Hausdorff, the inverse of $f$ is also continuous, so $f$ is a homeomorphism.
\end{proof}

Note that for $x \in \Xi$ we have $f(0,x) = x$, where $0 = (\mathrm{id},\mathrm{id})$ is the origin of $B$ with $x$ belonging to the interior of $X$. This establishes that $\Omega_r$ is locally a product $B \times \Xi$ about points of $\Xi$, with each chart $\phi_x$ given by $f^{-1}$ for $x \in \Xi$. By choosing appropriate values of $\epsilon$, the full regularity condition of the first part of Definition \ref{MatMan} (the extension of chart maps $\phi_x$ to $\hat{\phi}_x$) follows.  We can easily move these charts to other locations using the action of rigid motion. Indeed, take any $T' \in \Om_r$ and a rigid motion $g \in \R^d \rtimes \SO(d)$ so that $g(T') \in \Xi$. Such a rigid motion can be found, for example, by first translating a puncture of a tile in $T'$ over the origin, and then applying a rotation so as to orient this central tile as a prototile of $\mathcal{P}$. So we may consider the map
\[
f_g \colon B \times \Xi \to g^{-1}(X), \ f_g(b,Y) \coloneqq g^{-1}(b(Y)),
\]
which satisfies $f_g(0,g(T')) = T'$ with $g(T') \in \Xi$. Since rigid motions act as homeomorphisms on $\Om_r$, the above shows that there is a local product structure of $B \times \Xi$ at all points of $\Om_r$, where $B$ is homeomorphic to a $d+d(d-1)/2$-dimensional closed ball and $\Xi$ is a compact, separable, totally disconnected space and hence:

\begin{proposition}\label{OmrMM}
For a tiling $T$ with Euclidean FLC, $\Om_r$ is a matchbox manifold with Euclidean leaves of dimension $d+d(d-1)/2$.
\end{proposition}

We note that the situation is somewhat simpler than for the general case of a matchbox manifold since we can take each transverse model as the whole transversal $\Xi$, rather than a subspace depending on the chart. It is easily checked that transition maps between charts are isometries on the Euclidean coordinates, induced by multiplying by fixed rigid motions, which follows from the fact that all charts are related to the one about the canonical transversal by rigid motion. The plaques determined by the charts are given by small rigid motions of tilings, so the leaves, which are stitched together from these plaques, are given by orbits of tilings through rigid motion. As is always the case for matchbox manifolds, these leaves are precisely the path-components of $\Om_r$.

The above constructions extend those of the translational setting. When translational FLC is satisfied (which shall be assumed in the following chapters) $\mathcal{P}$ may be taken as the set of translation classes of tiles, with representatives taken with puncture over the origin. Then $\Xi$ is the standard canonical transversal. By omitting the rotational component in the above constructions, the space $\Om_t$ is equipped with charts making it a $d$-dimensional sub-matchbox manifold of $\Om_r$. The group $\SO(d)$ acts on $\Om_r$ transversally to the embedded subspace $\Om_t$, in the sense that sufficiently small rotates $g(\Om_t)$ are mutually disjoint. However, the subspaces $g(\Om_t)$ need not be disjoint for all distinct $g \in \SO(d)$. Non-trivial rotations of elements of $\Om_t$ can again belong to $\Om_t$. This can occur because of (discrete) rotational symmetries of tilings or, in general, because of elements of the \emph{point group}.


\section{Point Groups For FLC Tilings}\label{sec: rotations}

In this section we set out what it means for a tiling $T$ in $\R^d$ to have a point group, extending the usual notion from the periodic case. As in the previous section, given a rigid motion $\Phi$ we let $\Phi(T)$ denote the tiling given by applying $\Phi$ to each tile of $T$ (preserving their labels, if the tiles are labelled). Similarly, given a patch $P \subset T$ we can define the patch $\Phi(P) \subset \Phi(T)$.

\begin{definition}\label{def: acts} 
We say that a finite subgroup $G \leq \SO(d)$ \emph{acts on} $T$ \emph{by rotations} if, for every patch $P$ of $T$ and $g \in G$, we have that $g(P)$ is also a patch of $T$, up to translation.
\end{definition}
 
\begin{proposition} \label{prop: act}
If $G$ acts on $T$ by rotations then $G$ canonically acts on $\Omega_t$, by $g\cdot T' \coloneqq g(T')$.
\end{proposition}

\begin{proof} If $T' \in \Omega_t$ then every patch $P$ of $T'$ is a patch of $T$. Since $G$ acts on $T$ by rotations, $g(P)$ is also a patch of $T$. It follows that every patch of $g(T')$ is a patch of $T$, so $g(T') \in \Omega_t$. This is a continuous group action, from the definition of the topology of $\Om_t$. 
\end{proof}
 
We would like to define the point group of $T$ to be the maximal subgroup of $\SO(d)$ that acts. There will be such a maximal group, all groups that act are subgroups of the finite group of all rotational symmetries of the set of prototiles, but for general tilings this may not pass to a well defined group shared by all tilings in $\Om_t$: other elements of $\Om_t$ may have larger or smaller maximal groups acting on them.
 
\begin{example}\label{exp: bad}
Consider $T$ to be the tiling of $\R^2$ by unit squares, indexed by elements of $\Z^2$, where we colour squares $(x,y)$ black for $x = 0$ and white for $x \neq 0$. Then the group $\Z/2<\SO(2)$ consisting of the identity and rotation by $\pi$ acts on $T$ and its translates, but no larger subgroup of $\SO(2)$ does. Nevertheless, in $\Om_t$ there are periodic tilings consisting of only white tiles, and the group $\Z/4$ generated by rotation by $\pi/2$ acts on these.
\end{example}

Following \cite{Wal17rot} we may avoid issues related to this by demanding the point group only to be defined if a further condition is satisfied.
  
\begin{definition}\label{def: point group}
Say that $T$ \emph{has point group $G$} if $G$ acts on $T$ and there exists some $r > 0$ for which, whenever an $r$-patch $P$ as well as its rotate $g(P)$ belong to $T$ up to translation, for some $g \in \SO(d)$, then $g \in G$.
\end{definition}

In the example above, there is no point group defined because although $\Z/2$ is the largest rotation group that acts on $T$, there are patches $P$ of arbitrarily large radius for which $g(P)$ is in $T$ for $g$ rotation by $\pi/2$. It should be remarked, however, that it is not necessary for $\Om_t$ to contain only tilings with the same point group to ensure that the original has well defined point group. For example, tilings such as the next will not be excluded in what is to follow.

\begin{example} Consider a tiling $T$ of black and white unit squares of alternating concentric rings: a central tile is coloured black, its 8 neighbours are coloured white, their neighbours are coloured black, and so on. Then $\Z/4$ acts on $T$. However, there are tilings in $\Om_t$ of alternating black and white vertical lines of tiles, and similarly ones of horizontal lines. The point group for such tilings is $\Z/2$. More peculiar is that there are also limiting `corner tilings' in $\Om_t$ which do not have a well defined point group at all. Nonetheless, $T$ has well defined point group $\Z/4$ which still acts on $\Om_t$, permuting the $1$-periodic tilings and corner tilings between themselves.
\end{example}

Throughout the rest of the article we shall assume that our tilings have point groups; the following simple property of such tilings will be important in the next section in describing the topology of $\Omega_r$.

\begin{proposition} \label{prop: all rotations} Suppose $T$ has point group $G$. Let $T_1,T_2 \in \Omega_t$ with $g(T_1) = T_2$ for some $g \in \SO(d)$. Then $g \in G$.
\end{proposition}

\begin{proof}
Taking a patch $P$ at the origin of $T_1 \in \Om_t$, we have that $P \subset T$ must appear as a translate of a patch from $T$, and similarly $g(P) \subset T_2$ must appear as a translate in $T$. By taking $P$ sufficiently large, we see that $g \in G$.
\end{proof}

Propositions \ref{prop: act} and \ref{prop: all rotations} characterise the point group in terms of the action of rotation on $\Om_t$: suppose that $T$ has point group $G$ and that $g \in \SO(d)$. If $T' \in \Om_t$ and $g \in G$ then $g(T') \in \Om_t$. Conversely, if for some $T' \in \Om_t$ we have that $g(T') \in \Om_t$ too then $g \in G$. The same criteria may be used to define the point group for non-FLC tilings. For example, for the Pinwheel the tiling space $\Om_t$, taken as the completion of a translational orbit, is the full rotational hull of all pinwheel tilings, and one should take the point group to be the full rotation group $G = \SO(2)$.

Note that for a periodic tiling $\Om_t$ consists only of translates of $T$, so for $g \in \SO(d)$ we have that $g \in G$ if and only if $g(T)$ is a translate of $T$. Thus $G$ may be identified with the standard (orientation preserving) point group, i.e., as the quotient of (orientation preserving) symmetries of $T$ modulo translational symmetries. As addressed in Remark \ref{rem: reversing}, we restrict our attention to orientation preserving isometries in this article, but of course orientation reversing symmetries may easily be introduced in all of the constructions above.

\begin{remark} \label{rem: point group} In the case that the tiling has a point group, we can identify it as follows. Let $G_n$ be the subset of $\SO(d)$ of rotations $g$ such that whenever $P$ is an $n$-patch then $g(P)$ is also a patch of $T$, up to translation. It is easy to show, using FLC, that each $G_n$ is finite (indeed, note that the finite set of polyhedral prototiles have at most finite rotational symmetry). Moreover it is evident that each is closed under composition, so each $G_n$ is a subgroup of $\SO(d)$. They sit as a nested sequence 
\[
G_1\supset G_2\supset\cdots \supset G_n\supset\cdots
\]
since if any larger patch has a given rotation, then so must any smaller patch. Since each $G_n$ is finite the sequence is eventually constant, so there exists some $N$ for which $G_m = G_n$ for all $m,n\geqslant N$. It is then immediate from Definitions \ref{def: acts}, \ref{def: point group} that $T$ has point group $G_N$. \end{remark}

The next proposition shows that a repetitive tiling cannot be subject to the problems highlighted in Example \ref{exp: bad}. Since the tilings of main interest here are those which are highly structured, and certainly repetitive, the assumption that the tiling has a point group is therefore not restrictive in this setting.

\begin{proposition} \label{prop: repetitive => pg}
Suppose that $T$ is repetitive or, equivalently, that $(\Om_t,\R^d)$ is minimal. Then $T$ has a point group $G$. Moreover, every element of $\Om_r$ also has point group $G$.
\end{proposition}

\begin{proof}
Following Remark \ref{rem: point group}, we need to show that $T$ has point group $G_N$, where $N$ is chosen so that $G_m = G_N$ for $m \geqslant N$. By repetitivity, every $N$-patch appears in every $R$-patch for some $R>0$. We claim that if $P$ is an $R$-patch and $g(P)$ also appears in $T$ for $g \in \mathrm{SO}(d)$, then $g \in G_N$. Indeed, since every $N$-patch appears in $P$, every rotate by $g$ of patch of radius $N$ appears in $g(P)$ and hence in $T$. So by the definition of $G_N$ we have that $g \in G_N$, as required.

Every tiling $T' \in \Om_t$ has the same collection of finite patches, by repetitivity. Since the point group is defined purely in terms of these patches, every element of $\Om_t$ also has point group $G$. By Proposition \ref{prop: FLC => global rotate}, every element of $\Om_r$ is given by applying a global rotation to an element of $\Om_t$, which clearly does not change the point group.
\end{proof}

\subsection{Action of the point group on leaves}
As noted in Section \ref{sec: matchbox}, $\Om_t$ is a $d$-dimensional matchbox manifold, with leaves the translational orbits of tilings. Since $G$ acts on $\Om_t$ by rotations, it maps leaves diffeomorphically onto others. In fact, since a $d$-dimensional leaf consists precisely of the set of translates of some tiling, we have the following:

\begin{theorem}
Each $g \in G$ induces a bijection on the set of leaves of $\Om_t$. We have that $\ell = g(\ell)$ for a leaf $\ell$ if and only if each tiling $T' \in \ell$ has a global symmetry with rotational part given by $g$, that is, $T' = (\tau \circ g)(T')$, where $\tau$ is a translation $\tau = x \mapsto x+y$ for some $y \in \R^d$.
\end{theorem}

\begin{examples} The Penrose tilings have point group $\Z/10$, although we note that while there are two Penrose tilings with five-fold rotational symmetry about the origin, there are no individual tilings with ten-fold rotational symmetry. So the rotation $g$ by $2\pi/10$ is a bijection on the leaves of $\Om_t$ which sends leaves to distinct leaves, but $g^2$ sends precisely two leaves to themselves (each with a single fixed point: the unique tiling with $\Z/5$ rotational symmetry on each leaf). The chair tilings have point group $\Z/4$ and in this case there are examples with the full rotational symmetry group about the origin. The Ammann--Beenker tilings have point group $\Z/8$.
\end{examples}

We saw in Section \ref{sec: matchbox} that $\Om_r$ is a matchbox manifold, of dimension $d+d(d-1)/2$ whose leaves are the Euclidean orbits. In particular, since $\SO(d)$ is path-connected, for any $g \in \SO(d)$ and a leaf $\ell$, we have that $g \ell = \ell$. We have that $g(\Om_t) = h(\Om_t)$, as subspaces of $\Om_r$, if and only if $gh^{-1} \in G$. This observation implies a construction of $\Om_r$ from the action of $G$ on $\Om_t$, which we exploit in the following section.


\section{Useful fibrations}\label{sec: useful}

Our initial goal is to analyse the topology of the rotational hull $\Om_r$. To this end, in this section we introduce various fibre bundles which relate the topology of $\Om_r$ to the simpler translational hull $\Om_t$ and the action of the point group $G$ on it. In so doing we also introduce a further space, denoted $\Om_G$, the homotopy quotient of $\Om_t$ by the action of $G$.

Although our original motivation for $\Om_G$ was simply as a useful subsidiary space in studying $\Om_r$, its homotopy theory appears to be of fundamental importance, as we shall see in Section \ref{sec: homotopy}. We thus advise the reader to interpret the results of this section as relating the topologies of $\Om_r$ and $\Om_G$, our spaces of interest, to other spaces, namely the translational hull and various homogeneous spaces, which have been more widely studied.

The action of $\SO(d)$ on $\Omega_r$ is not free, nor is that of $G$ on $\Omega_t$. Indeed, tilings fixed by non-trivial rotational symmetries are of particular interest: it is precisely these which stop $\Om_r$ from being the simple product $\Om_t\times \SO(d)$. Topologically, however, it is easier to work with free actions, and a standard way of converting the $G$-action on $\Om_t$ into a free action is to add an additional factor upon which $G$ does act freely.

It turns out to be convenient in many cases to pass to the universal cover $\Sp(d)$ of $\SO(d)$, and the corresponding lift of $G$, which we shall denote $\GG$. For $d=2$, we take\footnote{It is common in the literature to denote the universal cover by $\widetilde{\Sp(2)}$ in dimension $d = 2$, and to take $\Sp(2) \to \SO(2)$ as the double fold cover $S^1 \xrightarrow{\times 2} S^1$. Since we always wish to consider the universal cover, for convenience we take the convention here that $\Sp(d)$ denotes the universal cover, even for $d = 2$.} $\Sp(2) \cong \R$ as the infinite cover of $S^1$, so $G \cong \Z/n$ for some $n \in \N$ and the elements of $\GG \cong \Z$ can be thought of as rotations by $2 \pi k /n$ for some $k \in \Z$, where rotates differing by a multiple of $2 \pi$ are distinguished in the lift $\GG$. For $d > 2$, $\GG$ sits as the upper extension in the diagram
\[
\begin{tikzcd}
\Z/2 \ar[r] \ar[d,equal] & \GG \ar[r,"q"] \ar[d] & G \ar[d]\\
\Z/2 \ar[r]  & \Sp(d) \ar[r]  & \SO(d)
\end{tikzcd}
\]
in which the right hand square is a pullback. Here, for $d>2$, $\Sp(d)$ is the universal double cover of $\SO(d)$, and we may think of its elements as simply rotations which are also imbued with an extra binary piece of information recording chirality (for example, a $2\pi$ rotation of $\R^3$ about an axis does not correspond to the identity, but a $4\pi$ rotation does). For $d=3$, we have that $\Sp(3) \cong S^3$ and the groups $\GG$ thus act freely on an odd dimensional sphere and so have periodic group cohomology \cite{TomZve08} (which need not be true for the original groups $G$, for example $G=A_5$). More generally, while $\SO(d)$ has non-trivial fundamental group, $\Sp(d)$ has $n^{\rm th}$ homotopy group $\pi_n(\Sp(d))=0$ for $n=1$ and $2$; this will give us computational advantage later.

The group $\GG$ still acts on $\Omega_t$, and $\Sp(d)$ acts on $\Omega_r$, in the obvious way via the quotient maps $q\colon \GG\to G$ and $\Sp(d)\to\SO(d)$, that is, by $g \cdot T \coloneqq q(g) \cdot T$: simply consider $g$ as a rotation, forgetting about the chirality information, then rotate in the usual way.

\subsection{First fibration for $\Om_r$}
Consider the action of $\GG$ on $\Omega_t \times \Sp(d)$ given by $g \cdot (T,s) \coloneqq (g \cdot T, sg^{-1})$. This is a free action, since it is free on the $\Sp(d)$ component. In fact, we have the following.

\begin{proposition} \label{prop: quotient} The quotient space $(\Omega_t \times \Sp(d))/\GG$ is homeomorphic to $\Omega_r$, so we have a principal $\GG$-bundle $\GG \to \Om_t \times \Sp(d) \to \Om_r$. \end{proposition}

\begin{proof} Denote the equivalence class in the quotient of any $(T,s) \in \Omega_t \times \Sp(d)$ by $[T,s]$. We let $f \colon (\Omega_t \times \Sp(d))/\GG \to \Omega_r$ be defined by $f([T,s]) \coloneqq s(T)$. The map is well defined since $g \cdot (T,s) = (g(T),sg^{-1})$ and $s(T) = (sg^{-1})(g(T))$.

The map is surjective, since every element of $\Omega_r$ is a rotate of some element of $\Omega_t$. To show that the map is injective, suppose that $f([T_1,s_1]) = f([T_2,s_2])$, so $s_2^{-1}s_1(T_1) = T_2$. Let $g \coloneqq s_2^{-1}s_1$. Since $g(T_1) = T_2 \in \Omega_t$, we must have that $g \in \GG$, since by Proposition \ref{prop: all rotations} the image of $g$ in $\SO(d)$ is in $G$. Moreover,
\[
g \cdot (T_1,s_1) = (s_2^{-1}s_1(T_1),s_1(s_2^{-1}s_1)^{-1}) = (T_2,s_2),
\]
so $[T_1,s_1] = [T_2,s_2]$. The quotient space is Hausdorff since $\GG$ is a proper group action; as $\Omega_r$ is compact $f$ is thus a homeomorphism. \end{proof}

\begin{remark}\label{rem: altfib}
In a similar way we have that $\Omega_r \cong (\Omega_t \times \SO(d)) / G$.
\end{remark}

\begin{corollary}\label{cor: firstfib}
There are fibrations
\begin{align*} \pushQED{\qed} 
& \Omega_t \times \Sp(d) \to \Omega_r \to B\GG \\
& \Omega_t \times \SO(d) \to \Omega_r \to B G.\qedhere
\end{align*}
\end{corollary}

Here $BG$, etc., denotes the {\em classifying space\/} of the group $G$. Recall that the set of principal $G$-bundles over a space $X$ are classified by the homotopy classes of maps $f\colon X\to BG$: given such a map $f$, the corresponding $G$-bundle is the pullback by $f$ of the universal $G$-bundle, which we write as $G\to EG\to BG$. As such, the universal bundle and the spaces $BG$ and $EG$ are only defined up to homotopy equivalence; indeed, given any contractible free $G$-space $E$, the quotient $E/G$ gives an example of a $BG$. The first (and similarly the second) fibration of the Corollary may be seen as the pullback of the universal bundle $\GG\to E\GG\to B\GG$ via the map $f \colon \Omega_r \to B\GG$ which classifies the $\GG$-bundle of Proposition \ref{prop: quotient} and which we shall also refer subsequently to as the {\em de-looping\/} of the bundle $\GG\to\Om_t\times \Sp(d) \to \Om_r$. Note that the action of $\pi_1(B\GG)=\GG$ on the fibre $\Omega_t \times \Sp(d)$ is (essentially by definition) precisely the action above.


\subsection{Second fibration for $\Om_r$}
Corollary \ref{cor: firstfib} gives the first of our two main fibrations for $\Om_r$. Its main application in this article will be in establishing Theorem \ref{thm: ratrot} on the rational cohomology of $\Om_r$, and deriving our second fibre bundle description of $\Om_r$. This second fibre bundle will be used for the deeper calculations of the integral cohomology groups in Section \ref{sec: planar} onwards. It has the advantage that the spaces involved are finite dimensional; by contrast, the first fibration has base space $B\GG$ or $BG$, which is typically cohomologically infinite dimensional over $\Z$-coefficients.

We consider the quotient space $\Sp(d)/\GG$ where $\GG$ acts on the right (so we identify $s$ and $sg$ for $s \in \Sp(d)$ and $g \in G$). We denote by $\zeta$ the quotient map $\Sp(d)\to \Sp(d)/\GG$.  We have a map $\theta \colon \Omega_r \to \Sp(d)/\GG$, which we call the \emph{orientation map}, defined as follows. Proposition \ref{prop: quotient} identifies $\Omega_r$ with $(\Omega_t \times \Sp(d)) / \GG$ and we define $\theta [T,g] = [g]$. Equivalently, every tiling $T$ of $\Omega_r$ is a rotate by some $g \in \Sp(d)$ of a tiling from $\Omega_t$, and $\theta(T) \coloneqq [g]$. The space $\Sp(d)/\GG \cong \SO(d)/G$ parametrises how patches are oriented with respect to those from the translational hull.

Recall the following construction of the fibre bundle associated to a principal $Q$-bundle and left $Q$-space (see \cite[Theorem 2.10]{MimTod91}):

\begin{lemma} \label{lem: diag fig} Let $Q$ be a topological group, $p \colon S \to S/Q$ a principal $Q$-bundle and $X$ a left $Q$-space. Then we have a fibre bundle
\[
X \longrightarrow \frac{X \times S}{Q} \mathop{\longrightarrow}^{p'} \frac{S}{Q}.
\]
where $p'[x,s] \coloneqq p(s)$ and $Q$ acts diagonally on $X \times S$. By choosing a base point of $S$, the fibre over $p(s)$ is the image of the canonical inclusion of $X$ into the first coordinate of $(X \times S)/Q$.
\end{lemma}

Recall from Corollary \ref{cor: firstfib} that we have a description for $\Om_r$ as the total space of a fibration with fibre $\Om_t \times \Sp(d)$ (or $\Om_t \times \SO(d)$) and base $B\GG$ (or $BG$, respectively). This comes from the de-looping of the principal $\GG$-bundle (or principal $G$-bundle) from Proposition \ref{prop: quotient}. Our second bundle, given by the orientation map, is derived from this by applying the above lemma. It expresses $\Om_r$ in terms of the finite dimensional spaces $\Om_t$ and $\Sp(d)/\GG$.

\begin{corollary} \label{cor: fibSOG} We have a fibre bundle
\[
\Om_t \to \Om_r \stackrel{\theta}{\to} \Sp(d)/\GG\,.
\]
The map $\theta$ is given by $\theta(g(T)) \coloneqq [g]$ for $T \in \Om_t$ and $g \in \Sp(d)$. \end{corollary}

\begin{proof} Define a left $\GG$-action on $\Sp(d)$ by $g \cdot s \coloneqq s g^{-1}$. Then the quotient $\Sp(d) \to \Sp(d)/\GG$ is a principal $\GG$-bundle, so by Lemma \ref{lem: diag fig} we have the fibre bundle
\[
\Om_t \to (\Om_t \times \Sp(d))/ \GG \xrightarrow{\theta} \Sp(d)/\GG
\]
where $\theta[x,s] \coloneqq [s]$ and the action of $\GG$ on the product is the one given before Proposition \ref{prop: quotient}. The element $[x,s]$ of the quotient is identified with $s(T) \in \Om_r$ under the homeomorphism constructed in the proof of Proposition \ref{prop: quotient} and the result follows. \end{proof}

\subsection{Approximant fibrations}
An important underlying perspective in much of the topological analysis of tiling spaces is that of Shape Theory, and to make use of this perspective it is necessary to have expansions of the spaces studied via \lq good\rq\ inverse limits of CW complexes, referred to as {\em approximants}, \lq good\rq\ here as usual merely meaning helpful for whatever issue is being studied. Following \cite{BDHS10} and in the spirit of \cite{AndPut98} we note that the translational hull may be written $\Om_t=\il K_n$ where the spaces $K_n,\ n\in\N$, are defined as $\Om_t/\!\!\sim_n$ for a set of equivalence relations $\sim_n$ on tilings given by 
\[
T_1\sim_n T_2 \quad \Longleftrightarrow\quad 
T_1\cap B_n(0)\equiv T_2\cap B_n(0)\,.
\]
Here $B_n(0)$ denotes the closed ball in $\R^d$ of radius $n$ centred at the origin: essentially, a tiling $T_1\in\Om_t$ is defined uniquely by the sequence of increasing patches $\{T_1\cap B_n(0)\}$. Equivalently, $\Om_t=\il \R^d/\!\!\approx_n$ where $\R^d$ is decorated with the single tiling $T$ and the relation $\approx_n$ is defined by
\[
x\approx_n y\quad \Longleftrightarrow\quad (T-x)\cap B_n(0)\equiv (T-y)\cap B_n(0)\,.
\]
The same equivalence relations apply also to the rotational tiling space $\Om_r$ and we define $J_n=\Om_r/\!\!\sim_n$, and similarly obtain the description $\Om_r=\il J_n$.

The fibrations above have analogues at the level of these approximants.

\begin{proposition}\label{prop: approx}
For sufficiently large $n$ there are fibrations
\[
K_n\times\Sp(d) \to J_n \to B\GG\qquad\quad\mbox{and}\qquad\quad K_n \to J_n \to \Sp(d)/\GG
\]
and similarly for $\Sp(d)$ and $\GG$ replaced by $\SO(d)$ and $G$. All maps involved commute in the obvious sense with the forgetful maps between approximants.
\end{proposition}

\begin{proof}
The proofs are the exact analogues of those for the corresponding fibrations of the complete tiling spaces. A point in $J_n$ is an $n$-patch in $\Om_r$, which is a rotation of an $n$-patch in $\Om_t$, unique in $\Sp(d)/\GG \cong \SO(d)/G$ so long as $n$ is sufficiently large, by the definition of the point group (given in Definition \ref{def: point group}).
\end{proof}

\subsection{The homotopy quotient of $\Om_r$ by rotations}
We finish this section with the introduction of a further space associated to a tiling $T$ and its rotation group. It will become important in the next section when we consider the homotopy groups of tiling spaces. Recall that for any space $X$ equipped with an action of a group $Q$, the \emph{Borel construction}, $X_Q$, or \emph{homotopy quotient} of $X$ is defined as $(X \times EQ)/Q$ where, as before, $EQ$ is any contractible free $Q$-space. As the space $EQ$ is only defined up to homotopy, so is $X_Q$. Note that if $Q$ is a subgroup of $S$, then any candidate for $ES$ is also an $EQ$. 

\begin{definition}
Given a tiling $T$ with point group $G$, define $\Om_G$ as the Borel space $(\Om_t\times EG)/G$ and $\Om_\GG$ as the Borel space $(\Om_t \times E\GG)/\GG$.
\end{definition}

We can alternatively express these spaces as homotopy quotients of the rotational hull:

\begin{proposition} \label{prop: OmG alternative}
We have a homotopy equivalence between $\Om_G$ and the homotopy quotient of $\Om_r$ by $\SO(d)$, given by $(\Om_r \times E\SO(d))/\SO(d)$. Similarly $\Om_\GG \simeq (\Om_r \times E\Sp(d))/\Sp(d)$.
\end{proposition}

\begin{proof}
By Proposition \ref{prop: quotient} we have a homeomorphism $\Om_r \cong (\Om_t \times \SO(d))/G$. With respect to this identification, the action by rotations of $\SO(d)$ on $\Om_r$ is given by $h \cdot [T,s] = [T,hs]$. So the Borel space of $\Om_r$ mod rotations is
\[
(\Om_r)_{\SO(d)} \cong \frac{\frac{\left(\Om_t \times \SO(d)\right)}{G} \times E}{\SO(d)},
\]
where $E = E \SO(d)$. This is the quotient of $(\Om_t \times \SO(d)) \times E$ by the equivalence relation $((T,s),e) \sim_a h \cdot (g \cdot(T,s),e) = ((gT,hsg^{-1}),he)$, where $g \in G$ and $h \in \SO(d)$. Consider instead the equivalence relation given by $((T,s),e) \sim_b ((gT,s'),ge)$, where $g \in G$ and $s' \in \SO(d)$ are arbitrary. Since $E G$ can be taken as $E = E \SO(d)$, this gives quotient $\Om_G$, since the middle coordinate is annihilated and the outer two are given the diagonal action by $G$ on $\Om_t \times E$.

These quotient spaces are homeomorphic, which may be seen as follows. We have a self-homeomorphism $f$ on $(\Om_t \times \SO(d)) \times E$ given by $((T,s),e) \mapsto ((T,s),s^{-1}e)$, with inverse $((T,s),e) \mapsto ((T,s),se)$. Now, $f$ identifies
\[
((T,s),e) \leftrightarrow ((T,s),e') \ , \ \ ((gT,s'),he) \leftrightarrow ((gT,s'),ge')
\]
where $g \in G$ and $h \in \SO(d)$ are arbitrary, $s' = hsg^{-1}$ and $e' = s^{-1}e$. The left-hand pairs are precisely those identified with the first equivalence relation, and the right-hand pairs are those identified with the second, so $f$ interleaves the equivalence relations and induces a homeomorphism of the quotients, as required. The proof for $\Om_\GG$ is analogous.
\end{proof}

\begin{remark}
Since we may take $\Om_G$ as the homotopy quotient of $\Om_r$ by the action of $\SO(d)$, this provides the appropriate analogue of this space when extending to the setting that $T$ has finite local complexity with respect to all rigid motions, rather than just translations (such as for $T$ a pinwheel tiling).
\end{remark}

The approximants $K_n$ of $\Om_t$ naturally provide approximants for $\Om_G$ and $\Om_\GG$. Indeed, $G$ and $\GG$ canonically act on each $K_n$ (a point of $K_n$ corresponds to an $n$-patch with prescribed origin, the point group acts by rotating such a patch about the origin). So we may take the associated Borel constructions $(K_n \times EG)/G$ and $(K_n \times E\GG)/\GG$. As $n$ varies these have obvious `forgetful maps' induced by those of $K_n$ and the identity on the second factors. It is easily verified that their inverse limits are homeomorphic to models for $\Om_G$ and $\Om_\GG$, respectively.

We summarise the main structures of this section in the following omnibus proposition. Again, note that our main interest is in the spaces $\Om_r$ and $\Om_G$, which the below shows may be studied via the action of the point group on $\Om_t$.

\begin{proposition} \label{prop: fib diagram} 
\begin{enumerate}
\item We have a commutative diagram
\[
\begin{tikzcd}
\GG \ar[r,equal] \ar[d] & \GG \ar[r] \ar[d] & \Om_t \ar[r,equal] \ar[d] & \Om_t \ar[d] \\
\Om_t \times E \Sp(d) \ar[d] & \ar[l]  \Om_t \times \Sp(d) \ar[r] \ar[d] & \Om_r \ar[r] \ar[d,"\theta"] & \Om_\GG \ar[d] \\
\Om_\GG & \ar[l] \Om_r \ar[r] & \frac{\Sp(d)}{\GG} \ar[r,"\xi"] & B\GG
\end{tikzcd}
\]
in which the vertical columns are fibre bundles, and principal $\GG$-bundles for the left-hand two. 
\item There is a similar diagram replacing occurrences of $\GG$ with $G$ and $\Sp(d)$ with $\SO(d)$.
\item Likewise, there are also such diagrams with the hulls replaced with their approximants, in which case the various structure maps between them commute with the maps of the diagram in the obvious sense.
\end{enumerate}
\end{proposition}

The $E \Sp(d)$ term of the left-hand fibre bundle can be ignored up to homotopy since it is contractible; the de-looping of this fibration corresponds to that on the right. The maps of this diagram are described in the course of its proof.

\begin{proof} By definition we may replace terms $\Om_\GG$ with $(\Om_t \times E \Sp(d))/\GG$ and, by Proposition \ref{prop: quotient}, $\Om_r$ with $(\Om_t \times \Sp(d))/\GG$. The first and second columns are then simply the principal $\GG$-bundles associated with the corresponding quotient maps. The third is the fibre bundle of Corollary \ref{cor: fibSOG} and the fourth is constructed analogously, using Lemma \ref{lem: diag fig}.

Following the above replacements the horizontal maps have the obvious description in each case. We take $T \in \Om_t$ as base point, and choose a base point $e \in \Sp(d)$. We have a map $\GG \to \Om_t$ defined by $g \mapsto g \cdot T$, and a map $\Sp(d) \to E \Sp(d)$ defined by $s \mapsto s^{-1} \cdot e$. The remaining horizontal maps are induced by these two and the diagram is easily seen to commute. The proofs with $G$ replacing $\GG$, or the hulls with their approximants, are essentially identical. \end{proof}

\begin{remark}
In earlier work \cite{BDHS10, Wal17rot} on rotational structures there is a yet further space $\Om_0$ considered, defined as $ \Omega_t / G$ or equivalently as $\Omega_r/\text{SO}(d)$. We do not discuss this space here, but we do note that $\Om_G$ is in some sense a homotopy analogue: whereas $\Om_0$ is the quotient of $\Om_t$ by the $G$-action, $\Om_G$ is the {\em homotopy quotient\/} of $\Om_t$ by $G$; homotopy quotients are more natural to consider in the context of homotopy invariants such as cohomology or homotopy groups.
\end{remark}


\section{Homotopy groups}\label{sec: homotopy}

In this section we use the fibrations introduced above to study and determine the relations between homotopical invariants of the spaces involved. We begin in the general setting, which covers both periodic and aperiodic examples. In Subsection \ref{phtpy} we specialise to the periodic setting, where it will transpire that the constructed invariants may be identified with the space groups of the given patterns. Then Subsection \ref{ahtpy} addresses our notion of {\em \asg}\  for aperiodic patterns, and demonstrates methods of computing this and related homotopical invariants. Finally, in Subsection \ref{casg} we relate our invariants to the crystallographers' aperiodic space group.

\subsection{Homotopy groups for aperiodic patterns} For aperiodic patterns the classical homotopy groups are not appropriate objects to consider, due to the pathological nature of the associated hulls. Instead one should utilise the perspective of Shape Theory, using the pro-homotopy groups $\ph_k(-)$ as replacements. We shall briefly digress to describe the bare-minimum framework and definition of these, see \cite{ClHu12, MarSeg82} for further details, in particular the latter for the full set-up and justification of the formal shape category.

For us, an inverse system of groups or {\em pro-group\/} is a diagram of groups and homomorphisms ({\em structure maps\/}) indexed over the natural numbers of the form
\[
G_0\longleftarrow G_1\longleftarrow G_2\longleftarrow \cdots \longleftarrow G_n\longleftarrow \cdots\,.
\]
It is useful to consider an individual group $G$ also as a pro-group by representing it as the inverse limit of copies of $G$ with the identity as each of the structure maps. A morphism between two pro-groups, $G_*$ and $H_*$, is a sequence of homomorphisms $g_i \colon G_i \to H_{j(i)}$, for some monotonically ascending sequence of integers $j(i)$ with $j(i)\to\infty$ as $i\to\infty$, commuting in the obvious way with the structure maps for $G_*$ and $H_*$. Such a morphism is {\em level-preserving\/} if $j(i)=i$ for all $i$; this is not required, but most of our examples will be of this form. Two pro-groups $G_*$ and $H_*$ are {\em pro-equivalent\/} if there are morphisms $g_*\colon G_*\to H_*$ and $h_*\colon H_*\to G_*$ whose composites are equal to the structure maps in $G_*$ and $H_*$ respectively.  

We will need to consider exact sequences of pro-groups. For all our purposes we can restrict to level preserving morphisms; in this case a sequence of morphisms $F_*\buildrel {f_*}\over\longrightarrow G_*\buildrel{g_*}\over \longrightarrow H_*$ is {\em exact at $G_*$} if each sequence $F_i\buildrel {f_i}\over\longrightarrow G_i\buildrel{g_i}\over \longrightarrow H_i$ is exact at $G_i$.  

Shape Theory considers representations of spaces $X$ in terms of {\em expansions}, diagrams involving (pointed) CW complexes $X_n$
\[
X_0\lla X_1\lla X_2\lla\cdots\lla X_n\lla\cdots \lla X
\]
satisfying a certain universal property with respect to $X$. An example of such an expansion would be when a space $X$ is given by the inverse limit $X=\il X_n$, but not all expansions need to be of this form.

The $k^{\rm th}$ {\em pro-homotopy group\/} of $X$, denoted $\ph_k(X,x)$, is then the pro-group given by the induced tower
\[
\pi_k(X_0,x_0)\lla \pi_k(X_1,x_1)\lla \pi_k(X_2,x_2)\lla\cdots\lla \pi_k(X_nx_n)\lla\cdots
\]
of homotopy groups, where $x$ is a base point of $X$, represented by  $(x_0,x_1,\ldots)$ in the tower. Choosing different inverse limit presentations can clearly change the specific components of the inverse system of homotopy groups, and as such one should consider it as only one of many possible presentations of $\ph_k(X,x)$ in the corresponding pro-category of groups, but all such presentations will be pro-equivalent.

More generally, Shape Theory considers a notion of {\em shape equivalence\/} of expansions. Loosely speaking, two spaces will be shape equivalent if there are morphisms of the expansions of the spaces which change the component CW complexes and structure maps only by homotopies.  Clearly two shape equivalent spaces will then share the same pro-homotopy groups, even though the resulting point set inverse limits of their expansions could be vastly different.

\begin{example}\label{stableeg}
The translational hull $\Om_t$ of the one-dimensional Fibonacci tilings may be constructed \cite{AndPut98, FoHuKe02} as an inverse limit involving a single space and structure map
\[
\begin{tikzcd}
Z & \ar[l,swap, "\gamma"] Z & \ar[l,swap, "\gamma"] Z & \ar[l,swap, "\gamma"] Z & \ar[l,swap, "\gamma"] \cdots
\end{tikzcd}
\]
The space $Z$ may be chosen to be homotopy equivalent to $W$, the one point union of two copies of the circle. Moreover, the structure map is itself a homotopy equivalence (though certainly not a homeomorphism). Thus there is a shape equivalence of $\Om_t$ to $W$, considered as having the constant expansion. I.e., there is a homotopy commutative diagram
\[
\begin{tikzcd}
Z \ar[d] & \ar[l,swap, "\gamma"] Z \ar[d] & \ar[l,swap, "\gamma"] Z \ar[d] & \ar[l,swap, "\gamma"] Z \ar[d] & \ar[l,swap, "\gamma"] \cdots \\
W & \ar[l,swap, "\mathrm{id}"] W & \ar[l,swap, "\mathrm{id}"] W & \ar[l,swap, "\mathrm{id}"] W & \ar[l,swap, "\mathrm{id}"] \cdots
\end{tikzcd}
\]
where the vertical arrows are homotopy equivalences.

Although $\Om_t$ is clearly far from being homeomorphic, or even homotopy equivalent to $W$, this diagram shows that their pro-homotopy groups are pro-isomorphic. In a similar manner, one may show that $\Om_t$ for any Sturmian tiling is shape-equivalent to $W$, although two such translational hulls are typically not homeomorphic.
\end{example}

The important aspect of this sort of observation for us is that it is at times convenient to change one of our tiling spaces to a shape equivalent, but much simpler one; this will not change the value of any shape invariants we apply, such as \v{C}ech cohomology, or pro-homotopy groups, but it can make the underlying topology much easier to manage. We shall see this working in practice later on, for example in Example \ref{exp: fibcube}.

\bigskip
We turn now to consider the pro-homotopy groups for the tiling spaces associated to our tiling $T$ of interest. For our application, we shall assume that  $T$ is used to point the hulls $\Omega_t$, $\Omega_r$ (and $\Om_G$, $\Om_\GG$, once an arbitrary base point of $E \Sp(d)$ is chosen), as well as all their approximants. It should be remarked that a change of base points may result in different pro-homotopy groups, related issues are considered in \cite{GelPro95}.

\begin{proposition} \label{prop: rot higher htpy Om_r} For $k>1$ we have isomorphisms $\ph_k(\Om_r,T) \cong \ph_k(\Om_t,T) \oplus \pi_k(\SO(d))$. \end{proposition}

\begin{proof} Consider the second column of Proposition \ref{prop: fib diagram} and de-loop to the fibration $\Om_t \times \Sp(d) \to \Om_r \to B\GG$. Following Proposition \ref{prop: approx} this may be expanded as a map of inverse systems
\[
\begin{tikzcd}[column sep=small]
\cdots & \ar[l] K_{n-1} \times \Sp(d) \ar[d] & \ar[l] K_n \times \Sp(d) \ar[d] & \ar[l] K_{n+1} \times \Sp(d) \ar[d] & \ar[l] \cdots & \Om_t \times \Sp(d) \ar[d] \\
\cdots & \ar[l] J_{n-1} \ar[d] & \ar[l] J_n \ar[d] & \ar[l] J_{n+1} \ar[d] & \ar[l] \cdots & \Om_r \ar[d] \\ 
\cdots & \ar[l,equal] B\GG & \ar[l,equal] B\GG & \ar[l,equal] B\GG  & \ar[l,equal] \cdots & B\GG\,. \\
\end{tikzcd}
\]
Applying $\pi_k(-)$, the proof follows from $\pi_k(B\GG) \cong 0$ and $\pi_k(\Sp(d)) \cong \pi_k(\SO(d))$ for $k>1$, the latter by the fact that by definition $\Sp(d)$ is the universal cover of $\SO(d)$. \end{proof}

\begin{proposition} \label{prop: rot higher htpys} For $k > 1$ we have isomorphisms $\ph_k(\Om_t,T) \cong \ph_k(\Om_G,T) \cong \ph_k(\Om_\GG,T)$. \end{proposition}

\begin{proof} The proof is similar to the above, where we use the right-hand fibrations of Proposition \ref{prop: fib diagram}. \end{proof}

The above two results shows that the higher pro-homotopy groups of the spaces $\Om_r$, $\Om_G$ and $\Om_\GG$ are functions only of those of $\Om_t$ and the Lie groups $\SO(d)$, the latter being closely related to the (unsolved problem of the) higher homotopy groups of spheres $\pi_n(S^m)$.

\begin{remark}\label{rem: hardhomotopy} In the periodic case $\Om_t$ is a $d$-torus and thus has no higher homotopy (see Theorem \ref{thm: per higher hom}). The higher pro-homotopy groups for aperiodic tilings are in general more complicated to describe. For an illustration, a dimension $d$, codimension 1 canonical projection tiling is shape equivalent to a punctured $d+1$ torus \cite{FoHuKe02}, $(\T^{d+1}\!-\!\{p\})$, and we may identify the pro-homotopy groups of the tiling space with the usual homotopy groups of the punctured torus. For $d\geqslant 2$, the inclusion $i\colon (\T^{d+1}\!-\!\{p\})\lra\T^{d+1}$ gives an isomorphism on $\pi_1$. However, if, like $\T^{d+1}$, the subspace $\T^{d+1}\!-\!\{p\}$ had all higher homotopy groups zero, the map $i$ would be a homotopy equivalence, by Whitehead's theorem. This is clearly not the case, for example by considering homology in degree $d+1$.
\end{remark}

To obtain results on the pro-fundamental group of the rotational hull we shall consider the third column of Proposition \ref{prop: fib diagram}, so we first determine the homotopy groups of $\Sp(d)/\GG$.

\begin{lemma} \label{lem: homog htpy} The map $\xi \colon \Sp(d)/ \GG \to B\GG$ of Proposition \ref{prop: fib diagram} induces an isomorphism in $\pi_1$. We have that $\pi_2(\Sp(d)/\GG) \cong 0$ and $\pi_k(\Sp(d)/\GG) \cong \pi_k(\Sp(d)) \cong \pi_k(\SO(d))$ for $k>1$. \end{lemma}

\begin{proof} We have a covering map $\Sp(d) \to \Sp(d)/\GG$ and, since $\Sp(d)$ is the universal cover of $\SO(d)$, we have that $\pi_k(\Sp(d)) \cong \pi_k(\SO(d))$ for $k > 1$. The final claim follows easily from these two facts.

There are canonical fibre bundles $\SO(d) \to \SO(d+1) \to S^d$, where an element $g \in \SO(d+1)$ is sent to $g(x_0) \in S^d$ for $x_0 \in S^d$ a chosen base point. A simple inductive argument then shows that $\pi_2(\SO(d)) \cong 0$. It follows from the above that $\pi_2(\Sp(d)/\GG) \cong 0$.

The map $\xi \colon \Sp(d)/\GG \to B\GG$ is the fibration classifying the principal $\GG$-bundle $\GG \to \Sp(d) \to \Sp(d)/\GG$. That $\xi_*$ induces an isomorphism in $\pi_1$  follows from the resulting long exact sequence in homotopy.  \end{proof}

\begin{theorem} \label{thm: rot fund} We have a diagram of pro-fundamental groups in which each row is a short exact sequence induced by a fibre bundle of Proposition \ref{prop: fib diagram}.
\[
\begin{tikzcd}
\ph_1(\Om_t,T) \ar[r] \ar[d,equal] & \ph_1(\Om_r,T)   \ar[r,"\theta"] \ar[d,"\cong"] & \GG \ar[d,equal] \\
\ph_1(\Om_t,T) \ar[r] \ar[d,equal] & \ph_1(\Om_\GG,T) \ar[r] \ar[d]                  & \GG \ar[d,"q"]     \\
\ph_1(\Om_t,T) \ar[r]              & \ph_1(\Om_G,T)   \ar[r]                         & G \,.
\end{tikzcd}
\]
\end{theorem}

\begin{proof} The quotient map $q \colon \GG \to G$ induces a map between the fibre bundles $\Om_t \to \Om_\GG \to B\GG$ and $\Om_t \to \Om_G \to BG$, with the induced map on $\pi_1$ from $B \GG$ to $BG$ being identified with $q$. Similarly to the proof of Proposition \ref{prop: rot higher htpy Om_r}, we may expand the spaces of this diagram in a natural way using the approximants of Proposition \ref{prop: approx}. Applying homotopy we obtain the bottom two rows of the above diagram and maps between them, using the fact that $\pi_2$ of $BG$ and $B\GG$ are trivial. We may stitch on the top three terms by considering the map of fibre bundles from the third to fourth column of Proposition \ref{prop: fib diagram}, where we use the fact that $\pi_2(\Sp(d)/\GG) \cong 0$ and $\xi \colon \Sp(d)/\GG \to B\GG$ induces an isomorphism in $\pi_1$, by Lemma \ref{lem: homog htpy}. \end{proof}


\subsection{Homotopy groups for periodic patterns}\label{phtpy} We now apply the above results in the periodic setting. In this case the hulls associated to our patterns are CW complexes -- in fact, $\Om_t$ is a $d$-torus and $\Om_r$ is a manifold of dimension $d + d(d-1)/2$. The pro-homotopy groups may be replaced with the classical homotopy groups\footnote{This follows from the exact sequences associated to taking inverse limits and their first derived functors $\varprojlim^1$, which will be discussed in the next subsection. When $X$ is a CW complex, and so the tower of approximants may be taken as the constant diagram $X \xleftarrow{\mathrm{id}} X \xleftarrow{\mathrm{id}} \cdots$, we may identify $\varprojlim \pi_k(X) \cong \pi_k(X)$ and $\varprojlim^1 \pi_k(X) \cong 0$.}: results of the previous section still apply, but occurrences of $\ph_k$ may be replaced with $\pi_k$. Since all spaces in this section are connected CW complexes, we omit base points from homotopy groups.

\begin{theorem} \label{thm: per higher hom} For $T$ a periodic pattern in $\R^d$ there are isomorphisms $\pi_k(\Om_r) \cong \pi_k(\SO(d))$ for $k>1$. \end{theorem}

\begin{proof} This follows directly from Proposition \ref{prop: rot higher htpy Om_r}. Indeed $\pi_k(\Om_t) \cong 0$ for $k>1$ since $\Om_t$ is a $d$-torus and thus has trivial higher homotopy. \end{proof}

\begin{remark} Note that in the case $d=2$, this shows that $\Omega_r$ is an Eilenberg--Mac Lane space, that is $\pi_k(\Omega_r) \cong 0$ for $k > 1$. Therefore, for $d=2$, to calculate the cohomology of $\Omega_r$ is to calculate the group cohomology of $\pi_1(\Om_r)$. This fundamental group will be determined in Corollary \ref{cor: per pi1 Om_r}. \end{remark}

Recall that $\Gamma$, the \emph{space group} or \emph{crystallographic group} of a periodic tiling $T$ is the group of all symmetries of $T$, a subgroup of the full Euclidean group $E(d) \cong \R^d \rtimes \mathrm{O}(d)$. Let us call the {\it positive space group} $\Gamma_+$ the subgroup of all orientation preserving symmetries of $T$; it has the \emph{standard extension} $\Z^d \to \Gamma_+ \to G$ as the group of translational symmetries of $T$ by the point group $G$.

\begin{theorem} \label{thm: periodic fund} We have the following commutative diagram in which each row is a short exact sequence. Each of the top three rows is induced by a fibre bundle of Proposition \ref{prop: fib diagram}. The bottom row is the standard extension of $\Gamma_+$.
\[
\begin{tikzcd}
\pi_1(\Om_t) \ar[r] \ar[d,equal] & \pi_1(\Om_r)   \ar[r,"\theta"] \ar[d,"\cong"] & \GG \ar[d,equal] \\
\pi_1(\Om_t) \ar[r] \ar[d,equal]         & \pi_1(\Om_\GG) \ar[r] \ar[d]                 & \GG \ar[d,"q"]     \\
\pi_1(\Om_t) \ar[r] \ar[d,"\cong"] & \pi_1(\Om_G)   \ar[r] \ar[d,"\cong"] & G \ar[d,equal] \\
\Z^d         \ar[r]        & \Gamma_+     \ar[r]        & G  \\
\end{tikzcd}
\]
\end{theorem}

\begin{proof} The upper three rows of the diagram follow from the results of the previous subsection, so we need only establish the bottom isomorphism of extensions.

The third row of the diagram is induced from the de-looping of the principal $G$-bundle $G \to \Om_t \times EG \to (\Om_t \times EG)/G$. Let $\Lambda \cong \Z^d$ denote the normal subgroup of $\Gamma_+$ consisting of translational symmetries of $\Gamma_+$. The map $x \mapsto T+x$ induces a $G$-equivariant homeomorphism $\R^d/\Lambda \to \Om_t$, where the $G$-action on the quotient is the one induced, as a subgroup, by $\Gamma_+$ acting on $\R^d$ by isometries. Let $\Gamma_+$ act on $EG$ by considering only the rotational part of a symmetry, so in particular $\Lambda$ acts trivially on $EG$. Hence, we have $G$-equivariant homeomorphisms $\Om_t \times EG \cong (\R^d/\Lambda) \times EG \cong (\R^d \times EG) / \Lambda$ and thus an isomorphism of principal $G$-bundles:
\[
\begin{tikzcd}
G \ar[r] \ar[d,equal] & \Om_t \times EG \ar[r] \ar[d,"\cong"]  & (\Om_t \times EG)/G \ar[d,"\cong"] \\
G \ar[r]        & (\R^d \times EG)/\Lambda \ar[r,"f"] & ((\R^d \times EG)/\Lambda)/G\,.
\end{tikzcd}
\]
Writing $E \coloneqq \R^d \times EG$ we have that $\Gamma_+$ acts freely on the contractible space $E$. The map $f$ is then just the induced quotient map $E/\Lambda \to E/\Gamma_+$, so the corresponding map of fundamental groups may be identified with the inclusion $\Lambda \hookrightarrow \Gamma_+$. \end{proof}

\begin{corollary} \label{cor: per pi1 Om_G} We have an isomorphism $\pi_1(\Om_G) \cong \Gamma_+$, and in fact an isomorphism between the standard extension of $\Gamma_+$ and the extension $\Z^d \to \pi_1(\Om_r) \to G$ realised by the fibration $\Om_t \to \Om_G \to BG$ of Proposition \ref{prop: fib diagram}. \end{corollary}

\begin{corollary} \label{cor: per pi1 Om_r} The fundamental group $\pi_1(\Om_r)$ is uniquely determined as the pullback under the quotient map $q \colon \GG \to G$ of the standard extension of~$\Gamma_+$ to the extension $\Z^d \to \pi_1(\Om_r) \to \GG$, realised by the fibration $\Om_t \to \Om_r \to \Sp(d)/\GG$ of Proposition \ref{prop: fib diagram}. \end{corollary}

\begin{remark}
Note that $\Gamma_+$ is the semi-direct product of $G$ with $\Z^d$ precisely if there is a point in $T$ with full point group symmetry $G$; in this case then $\pi_1(\Om_r)$ is likewise the semi-direct product of $\GG$ with $\Z^d$.
\end{remark}

\begin{remark} \label{rem: pi1 Om_r}
The above result shows that $\pi_1(\Om_r)$ is to $\Gamma_+$ as $\GG$ is to $G$. An element of $\pi_1(\Om_r)$ is represented by a continuously parametrised motion of tilings $(T_t)_{t \in [0,1]}$, with $T_0 = T_1 = T$. This determines an orientation preserving symmetry of $T$. Conversely, such a symmetry determines a based loop in $\Om_r$, unique up to homotopy when complemented with winding number information in $\Z$ for $d = 2$ or with chirality information in $\Z/2$ for $d>2$. Forgetting the explicit extension, we obtain the following corollary.
\end{remark}

\begin{corollary}\label{cor: sg cover}
For a periodic tiling $T$ in $\R^2$, the fundamental group is an extension
\[
\Z \to \pi_1(\Om_r) \to \Gamma_+\,.
\]
For a periodic tiling $T$ in $\R^d$ for $d\geqslant3$, the fundamental group $\pi_1(\Om_r)$ is an extension 
\[
\Z/2 \to \pi_1(\Om_r) \to \Gamma_+ \,.
\]
\end{corollary}


\begin{remark} \label{rem: periodic extension}
The proof of Theorem \ref{thm: periodic fund} shows that $\Omega_G \simeq B \Gamma_+$, so its higher homotopy groups are trivial. Notwithstanding the comments at the end of Section \ref{sec: rotations} for why we restrict to orientation preserving symmetries of our tilings, Corollary \ref{cor: per pi1 Om_G} has an immediate analogue for the whole space group $\Gamma$. Denote by $G_\pm$ the  extension of the point group $G$ to include reflections which act on $T$, in the appropriate sense. Then $\Gamma$ is an extension $0\to\Z^d\to\Gamma\to G_\pm\to1$ and by the same argument as for Theorem \ref{thm: periodic fund} we obtain $\pi_1(\Om_{G_\pm}) \cong \Gamma$.
\end{remark}

\subsection{Topological space groups}\label{ahtpy} The above results show that the classical space group of a periodic pattern may be described as the fundamental group of $\Om_G$, and the spin cover of the space group as the fundamental group of $\Om_\GG$ or $\Om_r$. As discussed, in the aperiodic setting the fundamental groups of these spaces are not appropriate objects to consider, but we do have the pro-fundamental groups.

Taking the inverse limit of an inverse system of groups defines a functor from the category of pro-groups to the category of groups. Applied to the tower representing $\ph_k(X,x)$, we obtain the inverse limit $\il \pi_k(X_n,x_n)$; these are shape invariants of $X$ \cite{MarSeg82} which we call the \emph{shape homotopy groups}. 

While it is generally more straightforward to describe the inverse limit of a pro-group than the pro-group itself, it contains potentially less information. Moreover, inverse limits in this setting are only half exact functors, and so have first derived functors $\il^1 \pi_k(X,x)$. In the case of tiling spaces these are the $L$-invariants $L_k$ of \cite{ClHu12}. Recall that while $\il^1 \pi_k(X,x)$ is an abelian group for $k>1$, $\il^1 \pi_1(X,x)$ is only a pointed set; moreover, it is either the one point set 1, or is uncountable \cite[Theorem 2]{mGMo92}.

The following definition extends the notion of the space group from the classical periodic setting, as seen from the previous subsection.

\begin{definition} \label{def: top sp gp}
Let $T$ be an aperiodic tiling in $\R^d$ with point group $G$ and extension $G_\pm$ that includes reflections where relevant. Define the {\em topological space group\/} $\Gamma$ of $T$ as $\il \pi_1(\Om_{G_\pm},T)$ and the {\em topological \spg\/} $\sg$ as the pro-group $\ph_1(\Om_{G_\pm},T)$. Similarly, define the {\em positive topological space group\/} $\Gamma_+$ as $\il \pi_1(\Om_G,T)$ and  the {\em positive topological \spg\/} $\psg$ as $\ph_1(\Om_G,T)$.
\end{definition}

As the pro-homotopy groups are homeomorphism invariants, we immediately have:

\begin{theorem} The groups and pro-groups $\sg$, $\Gamma$, $\psg$ and $\Gamma_+$ are S-MLD invariants of $T$. \end{theorem}

The following extends Corollary \ref{cor: sg cover} to the aperiodic setting:

\begin{proposition} For $d=2$ the shape fundamental group $\il \pi_1(\Om_r,T)$ is a $\Z$-cover
\[
\Z \to \il \pi_1(\Om_r,T) \to \il \pi_1(\Om_G,T)= \Gamma_+
\]
of the positive topological space group, and is a $\Z/2$-cover
\[
\Z/2 \to \il \pi_1(\Om_r,T) \to \il \pi_1(\Om_G,T)= \Gamma_+
\]
for $d > 2$.
\end{proposition}

\begin{proof} By Proposition \ref{prop: OmG alternative} we have that $\Om_G \cong (\Om_r \times E \SO(d))/\SO(d)$, with associated fibration $\SO(d) \to \Om_r \to \Om_G$. We thus have the following long exact sequence of pro-homotopy groups
\[
\cdots \to \ph_2(\Om_r,T) \to \ph_2(\Om_G,T) \to \pi_1(\SO(d)) \to \ph_1(\Om_r,T) \to \ph_1(\Om_G,T) \to 0 \,,
\]
since $\SO(d)$ is connected. As can be seen from the proof of Proposition \ref{prop: OmG alternative}, the map $\Om_r \to \Om_G$ here may be identified with $f \colon (\Om_t \times \SO(d))/G \to (\Om_t \times E\SO(d))/G$, which corresponds to the map of fibrations from the second to the first column of Proposition \ref{prop: fib diagram}. These de-loop to the fibrations $\Om_t \times \SO(d) \to \Om_r \to BG$ and $\Om_t \times E\SO(d) \to \Om_G \to BG$. Since $\pi_2(\SO(d)) \cong \pi_2(BG) \cong \pi_3(BG) \cong 0$, it follows that the map $\Om_r \to \Om_G$ induces an isomophism $\ph_2(\Om_r,T) \cong \ph_2(\Om_G,T)$, and thus the above long exact sequence gives the short exact sequence
\[
0 \to \pi_1(\SO(d)) \to \ph_1(\Om_r,T) \to \ph_1(\Om_G,T) \to 0.
\]
Since $\pi_1(\SO(d)) \cong \Z$ for $d=2$ and $\pi_1(\SO(d)) \cong \Z/2$ for $d>2$, the above gives the desired short exact sequence after passing to inverse limits. \end{proof}

Just as in the periodic case, we have that $\psg$ is an extension of the point group by the shape fundamental group of $\Om_t$, which in the periodic case corresponds to the lattice of translations. Thus in the pro-category, we have directly from Theorem \ref{thm: rot fund}:

\begin{corollary}
There is a short exact sequence of pro-groups
\[
\ph_1(\Om_t,T)\longrightarrow \psg \longrightarrow G\,
\]
\end{corollary}

However, upon passing to limits, we obtain the following $5$ term exact sequences; these follow by applying inverse limits to Propositions \ref{prop: rot higher htpy Om_r}, \ref{prop: rot higher htpys} and Theorem \ref{thm: rot fund}. In the below we omit the base points, taken as $T$ for each space.

\begin{corollary} \label{cor: 5 terms} For $k > 1$ we have isomorphisms
\[
\il \pi_k(\Om_t) \cong \il \pi_k(\Om_G) \cong \il \pi_k(\Om_\GG), \ \ \ \ \il \pi_k(\Om_r) \cong \il \pi_k(\Om_t) \oplus \pi_k(\SO(d))
\]
and the $L_k$ invariants of $\Om_t$, $\Om_G$, $\Om_\GG$ and $\Om_r$ all agree. In degree $1$ we have the following commutative diagram
\[
\begin{tikzcd}
0 \ar[r] & \il \pi_1(\Om_t) \ar[r] \ar[d,equal] & \il \pi_1(\Om_r)   \ar[r] \ar[d,"\cong"] & \GG \ar[d,equal] \ar[r,"\partial'"] & \il^1 \pi_1(\Om_t) \ar[r] \ar[d,equal] & \il^1 \pi_1(\Om_r)   \ar[r] \ar[d,"\cong"] & 1\\
0 \ar[r] & \il \pi_1(\Om_t) \ar[r] \ar[d,equal] & \il \pi_1(\Om_\GG) \ar[r] \ar[d] & \GG \ar[d,"q"] \ar[r,"\widetilde{\partial}"] & \il^1 \pi_1(\Om_t) \ar[r] \ar[d,equal] & \il^1 \pi_1(\Om_\GG) \ar[r] \ar[d] & 1   \\
0 \ar[r] & \il \pi_1(\Om_t) \ar[r] & \il \pi_1(\Om_G)   \ar[r] & G \ar[r,"\partial"] & \il^1 \pi_1(\Om_t) \ar[r]        & \il^1 \pi_1(\Om_G)   \ar[r]        & 1 \,.
\end{tikzcd}
\]
in which each five term row is exact.
\end{corollary}

Using similar techniques and considering the fibration $\Om_t \to \Om_{G_\pm} \to B G_\pm$, we also have an exact sequence
\begin{equation} \label{eq: 5 term space group}
0 \to \il \pi_1(\Om_t) \to \il \pi_1(\Om_{G_\pm}) \to G_\pm \xrightarrow{\partial^\pm} {\il}^1 \pi_1(\Om_t) \to {\il}^1 \pi_1(\Om_{G_\pm}) \to 1
\end{equation}
for a general aperiodic pattern. 

Thus the space group will be an extension of $\il \pi_1(\Om_t,T)$ by the point group precisely when the appropriate $\partial$ map in the above $5$-term exact sequences is the trivial map; similarly for orientation preserving and $\Om_r$ cover cases. Certainly this will hold when, for example, the $L_1$-invariant for $\Om_t$ is trivial.

\begin{example} \label{exp: RPT} By \cite{GaHuKe13} the translational hull $\Om_t$ of a rational projection method tiling in $\R^d$ is \emph{stable}, in the sense \cite{MarSeg82}, i.e., that it is shape equivalent to a finite CW complex. Thus the $L_1$-invariant of $\Om_t$ is trivial and we have extensions
\[
\il \pi_1(\Om_t,T) \to \il \pi_1(\Om_r,T) \to \GG \,\qquad\qquad \ \ \il \pi_1(\Om_t) \to \il \pi_1(\Om_{G_\pm}) \to G_\pm \,.
\]
\end{example}

\bigskip

When a translate of $T$ realises the whole point group as its group of symmetries, we may also deduce that the $\partial$ maps of Corollary \ref{cor: 5 terms} or Equation \ref{eq: 5 term space group} are trivial:

\begin{theorem} \label{thm: point group realised} Let $g \in \GG$ and suppose that a translate $T'$ of $T$ satisfies $g(T') = T'$. Let $\partial'$ be as in the statement of Corollary \ref{cor: 5 terms}. Then $\partial'(g) = 0$. Moreover, if there exists a translate of $T$ preserved under rotation by the entire group $\GG$, then $\partial' = 0$ and the extension $\il \pi_1(\Om_t,T) \to \il \pi_1(\Om_r,T) \to \GG$ is split. 

Analogous statements hold for each of the other $5$-term exact sequences of Corollary \ref{cor: 5 terms} and Equation \ref{eq: 5 term space group}. \end{theorem}

\begin{proof} Although it need not be true that a change of base points preserves the pro-homotopy groups of the hulls, it is easily verified that changing base point by a translate defines natural isomorphisms between them. So without loss of generality we take $T'=T$.

It suffices to prove the result for the case of the rotational hull $\Om_r$ since the vanishing of the other $\partial$ homomorphisms will follow. Similar arguments will apply for the analogous splitting results.

Take a path $\zeta'(g) \colon [0,1] \to \Sp(d)$ from the identity to $g$. Applying the quotient, this defines a loop $\zeta(g)$ in $\Sp(d)/\GG$ representing $g \in \GG \cong \pi_1(\Sp(d)/\GG)$. This also defines based loops $\zeta_n(g)$ in the approximants $J_n$, given by rotating in the approximant according to $\zeta'(g)$; these are loops since $T$ is assumed to be preserved under rotation by $g$. By construction these loops are preserved under application of the forgetful maps $J_n \to J_{n-1}$, so the sequence $(\zeta_n(g))$ defines an element of $\il \pi_1(\Om_r)$. Each has trivial translational part, in the sense that each may be written as the loop $t \mapsto [b_n,\zeta'(g)(t)]$, where we identify $J_n \cong (K_n \times \Sp(d))/\GG$ (for sufficiently large $n$) and $b_n$ is the base point of $K_n$ corresponding to $T$. By construction, the map $\theta_* \colon \il \pi_1(\Om_r) \to \GG$ at the approximant level sends each such loop to $g \in \GG$, so $g$ is in the image of $\theta_*$, and by exactness the kernel of $\partial'$.

In the case where $T$ has full symmetry group $\GG$, it is not hard to verify that the map $g \mapsto (\zeta_n(g))_n \in \il \pi_1(\Om_r)$ is a homomorphism, defining a splitting of the short exact sequence. The proofs for the other exact sequences are analogous. \end{proof}

It is tempting to conjecture that the above theorem still holds when the assumption that $T'$ is a translate of $T$ is replaced with $T'$ only being an `almost translate' of $T$, that is with $T' \in \Om_t$. However, we have heuristic reasons to believe that one should not expect for such a result to hold, although as yet do not have a full counter-example to the claim.

We finish this subsection with a couple of example computations of \asg s.

\begin{example}\label{exp: codim1} Consider canonical codimension 1 projection tilings in $\R^d$ whose tiles are formed by projecting the $d$-skeleton of the unit cubical tesselation with vertices $\Z^{d+1} < \R^{d+1}$ in a strip of the form $E + I^{d+1}$, where $E$ is a suitably irrationally positioned $d$-dimensional hyperplane in $\R^{d+1}$ and $I^{d+1}$ denotes the unit hypercube. The translational hulls of such patterns, as mentioned in Remark \ref{rem: hardhomotopy}, are shape equivalent to punctured $(d+1)$-tori. In particular the translational and rotational hulls are \emph{stable} in the sense of \cite{MarSeg82}, that is, shape equivalent to CW complexes, so their pro-homotopy groups are base point independent.

The projected $1$-skeleton has $1$-cells of $(d+1)$ types, one for each direction of edge in $I^{d+1}$. The full point group $G_\pm$, including orientation reversing symmetries, must permute the corresponding vectors, and their negatives, and it is not too hard to show that non-periodicity implies that in fact $G_\pm \cong \Z/2$, consisting only of the identity and $x \mapsto -x$. Since the latter is orientation preserving if and only if $d$ is even, we have that $G \cong \Z/2$ for $d$ even and $G \cong 0$ for $d$ odd. For $d$ even, following some further simple calculations from the observations of Remark \ref{rem: hardhomotopy}, one may show that the action of the generator of $G$ on $\il \pi_1(\Om_t,T) \cong \Z^{d+1}$ is given by $v \mapsto -v$. There are tilings of the hull with $\Z/2$ rotational symmetry, so by Theorem \ref{thm: point group realised}, for $d$ even the \asg\ $\Gamma=\Gamma_+=\il \pi_1(\Om_G,T)$ is the semi-direct product of $\Z/2$ and $\Z^{d+1}$. For $d$ odd this also describes $\Gamma$, but $\Gamma_+=\il\pi_1(\Om_t,T)$.

For $d$ odd $\Om_r\cong \Om_t\times \SO(d)$, so $\il\pi_1(\Om_r,T) \cong\Z^{d+1}\times \Z/2$. For $d$ even, the covering group $\GG$ is $\Z$ for $d=2$ and $\Z/4$ for $d\geqslant4$. In these cases $\il\pi_1(\Om_r,T)$ is the corresponding semi-direct product of $\GG$ and $\il\pi_1(\Om_t,T)$ with action of $\GG$ as given above via the quotient $\GG\to G$.
\end{example}

The following is an example where we can describe completely all the homotopy for a 3-dimensional aperiodic tiling.

\begin{example}\label{exp: fibcube}
Let $S$ be a dimension and codimension 1 canonical projection tiling, so a tiling of two tile types which occur according to a Sturmian sequence, equivalently a cutting sequence \cite{Ser85} of irrational slope. It is readily checked that the corresponding tiling space $\Om_S$ is closed under taking mirror images of sequences. Let $\FF$ be the aperiodic tiling on $\R^3$ given by cuboid tiles decorated as the product of three copies of the 1-dimensional  tiling $S$ and write $\Om_t(\FF)$, $\Om_G(\FF)$, etc., for its respective hulls. Then for all $n>1$,  $\il \pi_n(\Om_G(\FF),\FF) = \il \pi_n(\Om_t(\FF),\FF)=0$ and   $\il \pi_n(\Om_r(\FF),\FF) = \pi_n(\SO(3))$, while $\il \pi_1(\Om_G(\FF),\FF)$ is a semi-direct product 
\[
1 \lra F_2\oplus F_2\oplus F_2 \lra \il \pi_1(\Om_G(\FF),\FF) \lra G \lra 1\,.
\]
Here $F_2$ denotes the free group on two generators, and $G$ is the  group of rotational symmetries of the cube. That the point group of $\FF$ is the full rotation group of the cube, and that this is split (using Theorem \ref{thm: point group realised}) follows from the fact that $\Om_S$ carries the action of $\Z/2$ given by reflection and, as in the example above, the tiling space is stable and so these constructions are base point independent. The action of $G$ on $F_2\oplus F_2\oplus F_2$ is as follows. The group $G$ is generated by permutation matrices
\[
R\ =\,\left(\begin{array}{ccc}
0&1&0\\     -1&0&0\\      0&0&1
\end{array}\right)
\qquad\qquad
D\,=\,\left(\begin{array}{ccc}
0&1&0\\     0&0&1\\      1&0&0
\end{array}\right)\,.
\]
Let us write the generators of the three copies of $F_2$ as elements $a_1$ and $b_1$, $a_2$ and $b_2$ and $a_3$ and $b_3$. Then the elements $R$ and $D$ act 
\[
R\colon\ \left\{\begin{array}{lcl}
a_1&\mapsto &a_2\\
b_1&\mapsto &b_2\\
a_2&\mapsto &a_1^{-1}\\
b_2&\mapsto &b_1^{-1}\\
a_3&\mapsto &a_3\\
b_3&\mapsto &b_3
\end{array}\right.
\qquad\qquad
D\colon\ \left\{\begin{array}{lcl}
a_1&\mapsto &a_2\\
b_1&\mapsto &b_2\\
a_2&\mapsto &a_3\\
b_2&\mapsto &b_3\\
a_3&\mapsto &a_1\\
b_3&\mapsto &b_1
\end{array}\right.
\]

To see all this,  we use the fact that the 1-dimensional  tiling space $\Om_S$ is shape equivalent to $W$, the one point union of two circles (for example, see \cite{FoHuKe02} Chapter III), and hence $\Om_t(\FF)$ is shape equivalent to $W\times W\times W$. Thus $\il \pi_1(\Om_t(\FF),\FF)=F_2\oplus F_2\oplus F_2$ and $\il \pi_n(\Om_t(\FF),\FF)=0$ for $n>1$, by the homotopy of $W$ and that homotopy groups take cartesian products to direct sums. \end{example}

\subsection{Relation to the crystallographic aperiodic space group}\label{casg} In Crystallography there has long been a concept of {\em aperiodic space group\/} ($\!ASG$) that generalises the classical space group to the case of quasicrystals \cite{dWJaJa81}. Briefly, let us assume the quasicrystal is modelled on an irrational slice through the higher dimension lattice $\Lambda \cong \Z^k$, and denote by $\Lambda'$ the projection of $\Lambda$ to the physical space $\R^d$ of the tiling (isomorphic to $\Lambda$ as a group). Then the $ASG$ is the  extension of the full point group $G_\pm$ (as in Section \ref{sec: rotations}, and containing reflections as well as orientation preserving symmetries) by $\Lambda'$, using the natural action of $G_\pm$ on $\Lambda'$. Thus it is an extension 
\begin{equation} \label{eq: ASG}
\Lambda'\lra ASG\lra G_\pm\,.
\end{equation}
Here we elaborate the relation of the $ASG$ to the \asg\ $\Gamma=\il \pi_1(\Om_{G_\pm},T)$ and the pro-group $\sg =\ph_1(\Om_{G_\pm},T)$. 

Recall that the dynamical system $\Om_t$ with translation action by $\R^d$ has a maximal equicontinuous factor (MEF), here denoted by $\EE$, with  factor map denoted $\eta\colon \Om_t\to \EE$. In turn this induces a homomorphism in degree 1 cohomology $\eta^*\colon H^1(\EE;\Z)\to H^1(\Om_t;\Z)$. This homomorphism is examined in detail by Barge, Kellendonk and Schmieding \cite{BKS}. In particular, they prove that $\eta^*$ is injective and, in many important cases such as for {\em almost canonical projection tilings\/} (see \cite{GaHuKe13, Ju}), or when $H^1(\Om_t;\Z)$ is finitely generated, the image of $\eta^*$ is a direct summand. 

We now give a variant of the topological space group which, as we shall see, corresponds to the $ASG$. The action of $G_\pm$ naturally induces an action on the MEF $\EE$ (which is a $k$-torus for such projection tilings) and so we may define the Borel construction $\EE_{G_\pm}$ by $(\EE \times EG_\pm)/G_\pm$ as before.

\begin{definition} We let $\Gamma^{\mathrm{pro}}_{\EE} \coloneqq \ph_1({\EE}_{G_\pm})$ and $\Gamma_{\EE} \coloneqq \il\pi_1({\EE}_{G_\pm})$. \end{definition}

Analogously to the fourth column of Proposition \ref{prop: fib diagram} we have a fibre bundle
\begin{equation} \label{eq: ASG fb}
\EE \to \EE_{G_\pm} = (\EE \times EG_\pm)/G_\pm \to BG_\pm
\end{equation}
and the map $\eta$ induces a map between these fibrations, in particular a map $\sigma \colon \Om_{G_\pm} \to {\EE}_{G\pm}$. So we have a diagram of extensions
\begin{equation} \label{eq: TSG->TSG'}
\begin{tikzcd}
\ph_1(\Om_t,T) \ar[r] \ar[d,"\eta_*"] & \sg \ar[r] \ar[d,"\sigma_*"]       & G_\pm \ar[d,equal]\\
\ph_1(\EE,T)   \ar[r]                 & \Gamma^{\mathrm{pro}}_{\EE} \ar[r] & G_\pm \,.
\end{tikzcd}
\end{equation}

Before examining this diagram further, we would like to relate the crystallographers' aperiodic space group $ASG$ to $\Gamma_{\EE}$ within a broad setting in which the former is defined, the \emph{rational projection tilings} \cite[Section 4]{GaHuKe13}.

For the rational projection tilings the cohomology $H^1(\Om_t;\Z)$ is finitely generated and the lattice $\Lambda$ may be identified with the first {\em homology\/} $H_1(\EE)$ of the MEF \cite{FoHuKe02,GaHuKe13}. Its image under the projection to $\R^d$ associated to the pattern is $\Lambda'$, which can be identified with the translation module of the tiling, dual to the Fourier module $H^1(\EE;\Z)$. In fact, $\EE$ may equally be identified as the classifying space of $\Lambda$. Since $\pi_1(\EE)$ is abelian it is naturally identified with $H_1(\EE)$ under the Hurewicz homomorphism. So applying the long exact sequence in homotopy establishes the following:

\begin{theorem} Let $T$ be a rational projection tiling. Then there is a natural isomorphism between the extension of Equation \ref{eq: ASG} and the extension
\[
\il\pi_1(\EE) \to \Gamma_{\EE} \to G_\pm
\]
induced by the bottom row of Equation \ref{eq: TSG->TSG'}. In particular, for such a tiling $ASG \cong \Gamma_{\EE}$. \end{theorem}

An alternative argument for the above is to consider the $ASG$ as acting on the ambient Euclidean space of the lattice $\Lambda \leqslant \R^N$ and applying a similar line of reasoning to the proof of Theorem \ref{thm: periodic fund}, namely via the observation that we may take $\R^N \times EG_\pm$ as $E(ASG)$.

Having connected $ASG$ with $\Gamma_{\EE}$, we now want to compare the latter with the topological space group $\Gamma$, that is we consider the map $\sigma_*$ of Equation \ref{eq: TSG->TSG'} upon passing to the inverse limit. Let us now consider any repetitive tiling with stable tiling space. In Shape Theory, a space is termed {\em stable\/} if it can be written as an inverse limit of finite CW complexes
\begin{equation} \label{eq: stable}
X_1\buildrel \alpha_1\over\longleftarrow X_2\buildrel \alpha_2\over\longleftarrow X_3\buildrel \alpha_3\over\longleftarrow \cdots \longleftarrow X_n\buildrel \alpha_n\over\longleftarrow \cdots X
\end{equation}
whose structure maps $\alpha_n$ are homotopy equivalences, as in Example \ref{stableeg}. In this situation, there is no information lost passing from the pro-group $\ph_1(X,x)$ to its limit group $\il\pi_1(X,x)$. It follows (from stability and \cite{BKS}) that a tiling whose translational hull $\Om_t$ is stable has MEF $\EE$ a torus; moreover, its cohomology $H^1(\EE;\Z)$, which can be interpreted as the group of dynamical frequencies associated to the $\R^d$ action on $\Om_t$, is a finitely generated, free abelian group.

\begin{theorem}\label{stableASG}
Let $T$ be a tiling with $\Om_t$ stable and $H^1(\Om_t;\Z)$ finitely generated. Then the map $\sigma_* \colon \Gamma \to \Gamma_{\EE}$ induced by the corresponding map of Equation \ref{eq: TSG->TSG'} is a surjection.
\end{theorem}

\begin{proof}
By stability $\il^1\pi_1(\Om_t,T) = 0$ and, by \cite{BKS}, $\EE$ must be a torus so $\il^1\pi_1(\EE) = 0$ too. It follows that we may replace the pro-groups of Equation \ref{eq: TSG->TSG'} with their inverse limits:
\[
\begin{tikzcd}
\il\pi_1(\Om_t,T) \ar[r] \ar[d,"\eta_*"] & \Gamma \ar[r] \ar[d,"\sigma_*"]       & G_\pm \ar[d,equal]\\
\pi_1(\EE)   \ar[r]                     & \Gamma_{\EE} \ar[r] & G_\pm \,.
\end{tikzcd}
\]
So it is sufficient to show that $\eta_*$ is surjective.

Take a stable inverse limit expansion for $\Om_t$, as in Equation \ref{eq: stable}. By a theorem of Rogers \cite{Rogers} we can realise $\eta$, up to homotopy, as a map of expansions where, for example, we take the constant  expansion for $\EE$. So without loss of generality we may identify $\eta_* \colon \il\pi_1(\Om_t,T) \to \il\pi_1(\EE)$ with $f_* \colon \pi_1(X_1,x_1) \to \pi_1(\EE)$ for some suitable map $f \colon X_1 \to \EE$, where $X_1$ is the first CW approximant of the expansion. Because $H_1(\EE)$ is abelian, we may factor $f_*$ through homology
\[
\begin{tikzcd}
\il\pi_1(\Om_t,T) \cong \pi_1(X_1) \ar[r,two heads,"h"] & H_1(X_1) \ar[r,"f_*"] & H_1(\EE) \cong \pi_1(\EE) \,,
\end{tikzcd}
\]
by naturality of the Hurewicz homomorphism $h$. So $\eta_*$ is surjective if $H_1(X_1) \xrightarrow{f_*} H_1(\EE)$ is.

It follows from the universal coefficient theorem \cite{Spa66} that we may identify $H^1(\EE;\Z) \xrightarrow{\eta^*} H^1(\Om_t;\Z)$ with the dual
\[
\hom(H_1(\EE),\Z) \to \hom(H_1(X_1),\Z)
\]
of $f_*$. By \cite{BKS} $\eta^*$ is the injection of a direct summand. Since the groups involved are finitely generated abelian groups, it is easily verified that this can only happen if $f_*$ is surjective, as desired.
\end{proof}

Since rational projection method patterns have stable tiling spaces with finitely generated cohomology \cite{GaHuKe13}, we obtain the following:

\begin{corollary} Let $T$ be a rational projection tiling. Then there is a canonical surjection $\sigma_* \colon \Gamma \twoheadrightarrow ASG$. \end{corollary}

In the most general case, for example when $\EE$ is a solenoid and $H^1(\Om_t;\Z)$ is no longer finitely generated as a group, the translation module typically vanishes. One nevertheless still has the Fourier module, corresponding to $H^1(\EE;\Z)$, with its $G_\pm$-action. With similar constructions to before, the pro-group $\ph_1(\Om_t,T)$ will still determine 
\[
H^1(\Om_t;\Z) \cong \varinjlim \hom(\ph_1(\Om_t,T),\Z)
\]
and the extension of $G_\pm$ by $H^1(\EE;\Z)$ is a subgroup of the extension of $H^1(\Om_t;\Z)$ by $G_\pm$ by \cite{BKS}.

In conclusion, the pro-group $\sg$ is generally richer than the $ASG$ or its analogue to the degree that $\ph_1(\Om_t,T)$ is richer than $\pi_1(\EE)$. The \asg\ is of course generally distinct from the $ASG$, not least it being a non-abelian invariant on the translational component.

\begin{remark}\label{repvar}
In practice, the topological space group $\Gamma$ and its related objects are hard to compute in general: even the shape fundamental group $\il\pi_1(\Om_t,T)$ of the translational hull is unknown for most tilings in dimensions greater than 1. The observations of this section however suggest variants that still retain new information, but are more readily accessible to computation. We note two possibilities, as follows, but leave these constructions for further investigation elsewhere.

First, there is a homological analogue: as the point group $G_\pm$ acts on $\Om_t$ and hence $\il\pi_1(\Om_t)$, so it acts on the homology $\il H_1(\Om_t)$ and we may form the extension 
\[
\il H_1(\Om_t)\lra \Delta\lra G_\pm\,,
\]
a homological analogue for $\Gamma$. In the stable, finitely generated case, as in Theorem \ref{stableASG}, this gives an S-MLD invariant which again surjects onto the ASG, by the same argument as above. For the non-stable case the pro-group analogue, or the cohomological object
\[
\varinjlim  H^1(K_n;\Z)\lra \Delta^*\lra G_\pm
\]
would be suitable replacements.

The second variant, following the approach noted in \cite{Sad:PEC} and used to great effect by G\"ahler in his recent and ongoing classification of certain one dimensional tilings, is to examine the representation variety of $\Gamma$, that is 
\[
\mbox{Rep}(\Gamma;S)\ =\ \varinjlim \hom(\pi_1(J_n);S)
\]
for suitably chosen finite groups $S$.
\end{remark}


\section{Cohomology}\label{sec: cohomology}
We turn now to consider the \v{C}ech cohomology of the rotational hulls $\Om_r$. The fibrations of Section \ref{sec: useful}
\begin{align*}
\Om_t \times \Sp(d) & \to \Omega_r \to B\GG \\
\Om_t               & \to \Om_r    \to \Sp(d)/\GG
\end{align*}
give rise to two Serre type spectral sequences

\begin{align*}
H^n(B\GG; H^k(\Om_t \times \Sp(d);\Z)) \Longrightarrow & H^{n+k}(\Omega_r;\Z)\\
H^n(\Sp(d)/\GG; H^k(\Om_t;\Z))         \Longrightarrow & H^{n+k}(\Omega_r;\Z)\,.
\end{align*}
It is important to note that in both cases, the left hand groups are cohomology with {\em twisted coefficients}: there is an action of the fundamental group of the base space (in both cases the group $\GG$) on the cohomology of the fibre, and in part it is here that much of the subtlety of the computations arise.  Cohomology $H^*(B\GG;M)$ for any $\Z\GG$ module $M$ can of course be interpreted as group cohomology $H^*(\GG;M)$. See \cite{Bro82} for details of group cohomology and its calculation. 

Of course both these spectral sequences require knowledge of the more traditional tiling cohomology, $H^*(\Om_t;\Z)$, as well as the $\GG$-action on this object. This is often highly non-trivial, though there are now many techniques that can be effective for individual classes \cite{BDHS10, GaHuKe13, SadBook08}.

We will not dwell on the closely related cohomology $H^*(\Om_G;\Z)$, but using the corresponding fibration 
\[
\Om_t\to\Om_G\to BG
\]
the same techniques detailed below will yield calculations for these groups as well. Indeed, the arguments for navigating the spectral sequences for specific cases of $H^*(\Om_r;\Z)$ generally give what is needed for computation of the corresponding $H^*(\Om_G;\Z)$.

In the next section we examine the rational cohomology $H^*(\Om_r;\Q)$. It is derivable from the integral cohomology by the universal coefficient theorem \cite{Spa66}, which in this case says that we have a natural isomorphism
\[
H^n(X;\Q) \cong H^n(X;\Z) \otimes \Q.
\]
In particular, the rational cohomology loses torsion information, and as a result is a far less subtle invariant. When we turn to the integral cohomology calculations of our periodic and aperiodic cubical examples, we shall also compute the $\Z/2$-coefficient cohomology. By the universal coefficient theorem \cite{Spa66}, this is derivable from the integral cohomology via the following split exact sequence:
\[
0 \to H^n(X;\Z) \otimes \Z/2 \to H^n(X;\Z/2) \to \mathrm{Tor}_1^\Z(H^{n+1}(X;\Z);\Z/2) \to 0.
\]
The $\Z/2$-cohomology has the advantage that the cohomology groups are in fact fields over $\Z/2$, so the extension problems in the associated spectral sequences are trivial.


\subsection{Rational cohomology}\label{sec: rational}

The case of cohomology for $\Om_r$ with rational coefficients proves to be significantly more straightforward than the integral problem, for one reason because of the finiteness of $G$, and hence the vanishing of the group cohomology $H^n(G;\Q)$ for positive $n$. The following result applies to all dimensions $d\geqslant2$, and irrespective of whether the tiling is periodic or aperiodic.

\begin{theorem}\label{thm: ratrot}
Let $T$ be a tiling in $\R^d$ with point group $G$. Then its rational \v{C}ech cohomology satisfies
\[
H^n(\Om_r;\Q)\ \cong \bigoplus_{a+b=n} H^a(\Om_t;\Q)^G\otimes H^b(\SO(d);\Q)
\]
where $H^a(\Om_t;\Q)^G=H^a(\Om_t;\Q)^{\GG}$ denotes the $G$-invariant elements of $H^a(\Om_t;\Q)$.
\end{theorem}

\begin{proof}
We use the fibration of Corollary \ref{cor: firstfib}, which runs $\Omega_t \times \SO(d) \to \Omega_r \to B{G}$. The Serre spectral sequence is a spectral sequence computing $H^*(\Om_r;\Z)$ with $E_2$-page
\[
E_2^{ij}\ =\ H^i(G;H^j(\Om_t\times \SO(d);\Q))\,.
\]
The right hand object is group cohomology with twisted coefficients in the $\Q$-module 
\[
H^j(\Om_t\times \SO(d);\Q)\ =\ \bigoplus_{r+s=j}H^r(\Om_t;\Q)\otimes H^s(\SO(d);\Q)\,.
\]
The $G$-action on this is the diagonal one, but note that the action of any subgroup of $\SO(d)$ on $\SO(d)$ is homotopically trivial ($\SO(d)$ is path connected, so given any $x\in\SO(d)$ and $g\in G$, there is always a path back from $g\cdot x$ to $x$. Thus $H^s(\SO(d);\Q)^G=H^s(\SO(d);\Q)$).

However, for any coefficient $G$-module $M$, as $G$ is a finite group, $H^i(G;M)$ for $i>0$ is killed by multiplication by $|G|$. Thus if $M$ is a $\Q$-vector space, $H^i(G;M)=0$ for all $i>0$. So the $E_2$-page of our spectral sequence is  non-zero only in the column $i=0$. Moreover, as $H^0(G;M)=M^G$, the $G$-invariant elements of $M$, the column $E_2^{0,j}$ is the $G$-invariants 
\[
H^n(\Om_t\times \SO(d);\Q)^G\ =\ \bigoplus_{a+b=n} H^a(\Om_t;\Q)^G\otimes H^b(\SO(d);\Q)\,.
\]
As there can be no differentials or extension problems with a single column spectral sequence, this completes the proof.
\end{proof}

\begin{example} \label{exp: codim1 rat coh} The dimension $d$, aperiodic, canonical, codimension 1 tilings of Example \ref{exp: codim1} have translational tiling hulls shape equivalent to a once-punctured $(d+1)$-torus, which has rational cohomology
\[
H^a(\Om_t;\Q) \cong
	\begin{cases}
	\Q                               & \text{for } a=0; \\
	\Q^{\binom{d+1}{a}} & \text{for } 0 < a < d+1; \\
	0                                 & \text{otherwise.}
	\end{cases}
\]
For $d$ odd $G \cong 0$ so $H^a(\Om_t;\Q)^G = H^a(\Om_t;\Q)$ and $\Om_r \cong \Om_t \times \SO(d)$. For $d$ even $G \cong \Z/2$ and it is not hard to show that the non-trivial element of $G$ acts trivially on $H^a(\Om_t;\Q)$ for $a$ even, so that $H^a(\Om_t;\Q)^G = H^a(\Om_t;\Q)$; and by $x \mapsto -x$ for $a$ odd, so $H^a(\Om_t;\Q)^G = 0$. \end{example}

\begin{example}\label{exp: fibcubecoh}
Consider the decorated cube tiling $\FF$ of Example \ref{exp: fibcube}. Up to shape equivalence $\Om_t(\FF)$ is the product $\Pi$ of three copies of $W$, the one point union of two circles. The group $G$ of rotational symmetries is the group of orientation preserving symmetries of the cube, whose action on the 1-cells are as described in Example \ref{exp: fibcube}. 

Let us name the 1-cells in the three copies of $W$ by $\{x_{ij}\,;\,1\leqslant i\leqslant 3,\ j=1,2\}$. These are the 1-cells of $\Pi$; the 2-cells are those products $x_{i_1j_1}x_{i_2j_2}$ where $i_1\not=i_2$; there are twelve of these. There are eight 3-cells, the eight products $x_{1j_1}x_{2j_2}x_{3j_3}$. As a cell complex, the boundaries all vanish.  We get a description of the integral cohomology $H^*(\Om_t(\FF);\Z)$: in dimensions 0, 1, 2 and 3 it is respectively $\Z$, $\Z^6$, $\Z^{12}$ and $\Z^8$, and is zero in higher dimensions. The cells as described may be taken as generators of these groups. The rational cohomology $H^*(\Om_t(\FF);\Q)$ is obtained by tensoring with $\Q$.

To compute $H^*(\Om_r(\FF);\Q)$ we need to determine the $G$-invariants $H^*(\Om_t(\FF);\Q)^G$. For each 1 or 2-cell there is a rotation that acts as the involution $x\mapsto-x$, and so there can be no $G$-invariant elements of $H^1(\Om_t(\FF);\Q)$ and $H^2(\Om_t(\FF);\Q)$. The degree zero group $H^0(\Om_t(\FF);\Q)$ is necessarily $G$-invariant, while in $H^3$ the  $G$-invariant  submodule is of rank 4, spanned by the set
\[
\left\{\begin{array}{ll}
x_{11}x_{21}x_{31}, \qquad &x_{11}x_{21}x_{32}+x_{11}x_{22}x_{31}+x_{12}x_{21}x_{31}\\
x_{12}x_{22}x_{32}, \qquad &x_{11}x_{22}x_{32}+x_{12}x_{22}x_{31}+x_{12}x_{21}x_{32}
\end{array}\right\}\,.
\]
Finally, note that $H^n(\SO(3);\Q) \cong H^n(S^3;\Q) \cong \Q$ for $n=0,3$ and is 0 otherwise, so by Theorem \ref{thm: ratrot} we obtain
\[
H^n(\Om_r(\FF);\Q)\ =\ \left\{\begin{array}{ll}
\Q&\mbox{for }n=0\\
\Q^5&\mbox{for }n=3\\
\Q^4&\mbox{for }n=6\\
0&\mbox{otherwise. }
\end{array}\right.
\]
\end{example}

\subsection{Top degree rational cohomology} \label{sec: top degree}
A similar argument to the proof of Theorem \ref{thm: ratrot}, using the fibration $\Om_t\to\Om_G\to BG$, shows that $H^r(\Om_G;\Q) \cong H^r(\Om_t;\Q)^G$. Likewise, for the space $\Om_0 \coloneqq \Om_t/G$ (equivalently, the quotient $\Om_r / \SO(d)$) it is not hard also to show that $H^r(\Om_0;\Q) \cong H^r(\Om_t;\Q)^G$ (c.f., \cite[Theorem 7]{BDHS10}). Theorem \ref{thm: ratrot} thus says that $\Om_r$ and $\Om_0 \times \SO(d)$ have isomorphic rational cohomology, generalising \cite[Theorem 8]{BDHS10}.

Let $D = d+d(d-1)/2$ (the top non-trivial cohomological degree of $\Om_r$). As a consequence of Theorem \ref{thm: ratrot}, we may relate the top degree $\Q$ or $\R$-coefficient cohomology of $\Om_r$ to the top degree cohomologies of $\Om_0$ and $\Om_t$.

\begin{corollary} \label{cor: top dim cohomology}
We have natural identifications
\[
H^D(\Om_r;\Q) \cong H^d(\Om_t;\Q)^G \cong H^d(\Om_0;\Q).
\]
\end{corollary}

\begin{proof}
Since $\SO(d)$ is a Lie group, its top degree cohomology group is $H^{d(d-1)/2}(\SO(d);\Q) \cong \Q$. By Theorem \ref{thm: ratrot}, we have that $H^d(\Om_t;\Q)^G \cong H^D(\Om_r;\Q)$. For the isomorphism $H^d(\Om_t;\Q)^G \cong H^d(\Om_0;\Q)$, see for example \cite[Proposition 3.12]{Wal17pe}.
\end{proof}

The top degree cohomology of $\Om_t$ has an important trace function for sufficiently regular tilings (those with \emph{uniform patch frequencies}, see \cite{BaGr}). Under the pattern-equivariant (PE) formalism (see \cite{Kel03, KelPut06, Sad:PEC}) representatives $\psi$ of classes in $H^d(\Om_t;\Q)$ may be viewed as PE cochains. In top cohomological degree $d$ for $\Om_t$, we may average the value of a $\Q$-valued PE $d$-cochain $\varphi$ on $d$-cells intersecting any given $r$-ball (by dividing the sum by the volume of an $r$-ball) which converges to some $\tau(\psi) \in \R$ as $r \to \infty$, assuming that $T$ has uniform patch frequencies. The value $\tau(\psi)$ does not depend on the representative of the cohomology class of $\psi$ taken, and $\tau$ induces a well-defined homomorphism
\[
\tau \colon H^d(\Om_t;\Q) \to \R
\]
called the \emph{trace}. By Corollary \ref{cor: top dim cohomology} this also defines a trace on the top-degree rational cohomology of $\Om_r$. In the translational setting, this trace function is central in Bellissard's Gap Labelling Theorem \cite{BBG06}, regarding the spectral gaps of quasiperiodic potentials. Arguably, comparing patches up to rigid motion rather than just translations may be more natural for certain applications, so it may be of interest to investigate the tracial theory from this perspective further. It is possible that this natural trace on $H^D(\Om_r;\Q)$ may be extended to lower degrees, which is done in the translational setting using the Ruelle--Sullivan current \cite{KelPut06}.


\subsection{Integral cohomology of the planar tilings}\label{sec: planar}
Integral \v{C}ech cohomology should be expected to be a considerably more subtle invariant, as it will see the higher cohomology of the finite group $G$ involved. In this section we restrict to dimension $d=2$, for which reasonably complete answers can frequently be given; this is perhaps not least because the only finite subgroups of $\SO(2)=S^1$ are the cyclic ones $\Z/n$. The following, using the techniques of this paper, recover the final results of the second author's work \cite{Wal17rot}. It applies equally whether the tiling is periodic or not.

\begin{theorem} {\em \cite[Theorem 4.1]{Wal17rot}}
Suppose $T$ is a tiling in $\R^2$ with point group $G$. Then the \v{C}ech cohomology of $\Om_r$ in degree $n$ is an extension
\[
0\lra H^{n-1}(\Om_t;\Z)_G\lra H^n(\Om_r;\Z)\lra H^n(\Om_t;\Z)^G\lra0
\]
where, for a $G$-module $M$, we use the usual notation of $M^G$ to mean the $G$-invariant elements of $M$, and we write $M_G$ to mean the $G$-coinvariants of $M$, the quotient of $M$ by the submodule generated by elements of the form $m-gm$ as $m$ and $g$ run over $M$ and $G$ respectively.
\end{theorem}

\begin{proof}
As $\SO(2)$ is just the circle $S^1$, the quotient $\SO(2)/G$ is also a copy of $S^1$. The fibration of the orientation map, Proposition \ref{cor: fibSOG}, is thus
\[
\Om_t \to \Om_r \to S^1 = B\Z
\]
and the Serre spectral sequence of this has just two non-zero columns, namely, writing in terms of group cohomologies, $E_2^{0,n}=H^0(\Z;H^n(\Om_t;\Z))$ and $E_2^{1,n}=H^1(\Z;H^n(\Om_t;\Z))$. The $\Z$-action on $H^*(\Om_t;\Z)$ given by the holonomy of this fibration is given by the natural $G$-action and the quotient $\Z\to G$. As any $0^{\rm th}$ group cohomology $H^0(G;M)$ is the $G$-invariants of $M$, the first column is identified as $H^*(\Om_t;\Z)^G$. As $\Z$ is a Poincar\'e duality group, we have an isomorphism of the second column
\[
E_2^{1,n}\ =\ H^1(\Z;H^n(\Om_t;\Z))\ =\ H_0(\Z;H^n(\Om_t;\Z))\ =\ H^n(\Om_t;\Z)_G
\]
since any $0^{\rm th}$ group {\em homology\/} $H_0(G;M)$ is the group of $G$-coinvariants of $M$. 
There is no room for any differentials in this spectral sequence, and the result follows.
\end{proof}

\begin{remark}
It is not clear that all these extensions will always split, since there may be torsion in $H^2(\Om_t;\Z)$: for example in the T\"ubingen Triangle Tiling \cite{Gah02}. The remaining extensions do not have this problem and so in general we have
\[
H^n(\Om_r;\Z) =
\begin{cases}
\Z                                                                     		 & \mbox{for }n=0, \\
\Z\oplus H^1(\Om_t;\Z)^G                                               		 & \mbox{for }n=1, \\
\mbox{Extension} \quad H^1(\Om_t;\Z)_G \to H^2(\Om_r;\Z) \to H^2(\Om_t;\Z)^G & \mbox{for }n=2, \\
H^2(\Om_t;\Z)_G                                                              & \mbox{for }n=3.
\end{cases}
\]
\end{remark}

\begin{example} Consider an aperiodic, canonical codimension $1$ projection tiling of $\R^2$. As in Example \ref{exp: codim1 rat coh}, $H^n(\Om_t;\Z) \cong \Z$, $\Z^3$, $\Z^3$ for $n = 0$, $1$, $2$, respectively, and is trivial otherwise. The action of the point group $G = \Z/2$ on $H^n(\Om_t;\Z)$ is the trivial action in degrees $n=0,2$ and $x \mapsto -x$ in degree $1$. Then $H^n(\Om_t;\Z) \cong H^n(\Om_t;\Z)^G \cong H^n(\Om_t;\Z)_G$ for $n=0$, $2$; $H^1(\Om_t;\Z)^G \cong 0$ and $H^1(\Om_t;\Z)_G \cong (\Z/2)^3$. Hence
\[
H^n(\Om_r;\Z) =
\begin{cases}
\Z                   & \mbox{for }n=0 \\
\Z                   & \mbox{for }n=1 \\
(\Z/2)^3 \oplus \Z^3 & \mbox{for }n=2 \\
\Z^3                 & \mbox{for }n=3.
\end{cases}
\]
\end{example}


\subsection{Dimension 3: Configurations of the cubic lattice}\label{sec: cube}

We consider here the cohomology of the space $\Or$ for the case of the periodic cubical tiling in $\R^3$. A point in this space corresponds to a placement of a (unit) cubical tessellation of $\R^3$ at some specific position relative to the origin and at some specific orientation relative to the coordinate axes. Thus in this case $\Om_r$ is a 6-manifold: the translational hull $\Om_t$ (for example, the subspace of these tilings with cube sides parallel to the axes) is a 3-torus, and there are a further 3 degrees of rotational freedom.

The computation of $H^*(\Or;\Z)$  runs as follows. The point group $G$ is the group of symmetries of the cube, and the index 2 covering group $\GG$ of $G$ acts freely on $S^3$ with a fundamental domain $1/48^{\rm th}$ of the sphere. We compute using the fibration $\T^3=\Ot\lra \Or\lra S^3/\GG$ which gives our main spectral sequence
\[
H^*(S^3/\GG;H^*(\T^3;\Z))\Longrightarrow H^*(\Or;\Z)\,.
\]
It is important to note that the action of $\pi_1(S^3/\GG)=\GG$ on $H^*(\T^3;\Z)$ is not trivial, and thus the left hand term (the $E_2$-page of the spectral sequence) is `cohomology with twisted coefficients'. Specifically, it is non-trivial on $H^n(\T^3;\Z)$ precisely for $n=1$ and 2, where it can be read off directly from the natural action of $G$ on the cubical lattice.

It will turn out that this spectral sequence has 2-torsion, and in order to solve the resulting extension problems we compute in parallel the analogous spectral sequence with $\F$, as opposed to $\Z$, coefficients.

In order to compute $H^*(S^3/\GG;H^*(\T^3;R))$ (where $R=\Z$ or $\F$), we first compute the group cohomology $H^*(B\GG;H^*(\T^3;R))$. By virtue of $\GG$ being a discrete subgroup of $S^3$, this group cohomology is 4-periodic and it can be computed using an efficient resolution of $\GG$ as described in \cite{TomZve08}. We obtain the results as shown in Figure \ref{fig: HBGP}.

\begin{figure}
\begin{tikzpicture}
  \matrix (m) [matrix of math nodes,
    nodes in empty cells,nodes={minimum width=5ex,
    minimum height=5ex,outer sep=-5pt},
    column sep=-1ex,row sep=-2ex]{
                &      &     &     &   & \\
          k     &&&&&&&\\
          3     &  \Z &  0 & 2  & 0 & 48 & \cdots &\\
          2     &  0 & 2 &  4  & 2  & 0 & \cdots&\\
          1     &  0 & 2 &  4  & 2  & 0 & \cdots&\\
          0      &  \Z &  0 & 2  & 0 &  48 & \cdots& \\
    \quad\strut &   0  &  1  &  2  & 3  & 4& n \strut \\};
  \draw[thick] (m-1-1.east) -- (m-7-1.east) ;
\draw[thick] (m-7-1.north) -- (m-7-7.north) ;
\end{tikzpicture}
\begin{tikzpicture}
  \matrix (m) [matrix of math nodes,
    nodes in empty cells,nodes={minimum width=5ex,
    minimum height=5ex,outer sep=-5pt},
    column sep=-1ex,row sep=-2ex]{
                &      &     &     &   & &  & \\
          k     &&&&&&&\\
          3     &  \F &  \F & \F  & \F & \F & \cdots&\\
          2     &  \F &  \F^2 & \F^2  & \F & \F & \cdots&\\
          1     &  \F &  \F^2 & \F^2  & \F & \F & \cdots&\\
          0     &  \F &  \F & \F  & \F & \F & \cdots&\\
    \quad\strut &   0  &  1  &  2  & 3  &4& n \strut \\};
  \draw[thick] (m-1-1.east) -- (m-7-1.east) ;
\draw[thick] (m-7-1.north) -- (m-7-7.north) ;
\end{tikzpicture}
\caption{Tables of the group cohomologies $H^n(B\GG;H^k(\T^3;\Z))$ (left) and $H^n(B\GG;H^k(\T^3;\F))$ (right). Recall that $H^k(\T^3;R)=0$ for $k>3$. In the notation  2 denotes the group $\Z/2$, and $4$ denotes $\Z/4$.  All these are 4-periodic, in the sense that  $H^n(B\GG;H^k(T^3;R))=H^{n+4}(B\GG;H^k(\T^3;R))$ for all $n\geqslant 1$.}
\label{fig: HBGP}
\end{figure}

The groups $H^*(S^3/\GG;H^*(\T^3;R))$ may be easily deduced from these calculations using the following lemma.

\begin{lemma} Let $Q$ be a finite subgroup of $S^3$ and $M$ a $Q$-module. Then $H^n(BQ;M) \cong H^n(S^3/Q;M)$ for $n=0,1,2$. We have that $H^n(S^3/Q;M) \cong 0$ for $n>3$, and for $n=3$ we have a long exact sequence:
\[
0 \to H^3(BQ;M) \to H^3(S^3/Q;M) \to H^0(BQ;M) \to H^4(BQ;M) \to 0.
\]
\end{lemma}

\begin{proof}
As $S^3/Q$ is a 3-manifold, necessarily its cohomology vanishes in dimensions more than 3. For low dimensions we note that $Q$ is a 4-periodic group and so it has a resolution that is 4-periodic. A model for $BQ$ can be taken with the corresponding cell structure, and $S^3/Q$ as its 3-skeleton. Denoting by $C^n$ the cochain group of $BQ$ in dimension $n$ with $M$ coefficients, we get a short exact sequence of cochain complexes
\[
0\to C^{*-4} \buildrel u\over\longrightarrow C^*\buildrel i\over\longrightarrow C^*(S^3/\GG)\to 0\,.
\]
Here $C^*(S^3/Q)$ is the cochain complex of $S^3/Q$ with $M$ coefficients, equal to $C^*$ for $*\leqslant 3$ and 0 for $*>3$, and $u$ is the periodicity operator, identifying $C^n$ with $C^{n+4}$ for $n\geqslant0$. The homomorphism $i$ may be identified as that induced by the inclusion $S^3/Q\to BQ$. The lemma now follows from the resulting long exact sequence: for $n<3$ we have 
\[
0\to H^n(BQ;M)\buildrel {i^*}\over\longrightarrow H^n(S^3/Q;M)\to 0
\]
and finally for $n=3$ the exact sequence as stated in the lemma.
\end{proof}

From this we deduce the $E_2$-pages of the spectral sequences $H^*(S^3/\GG;H^*(\T^3;R))\Longrightarrow H^*(\Or;R)$, as shown in Figure \ref{fig: E2 periodic}.

\begin{figure}
\begin{tikzpicture}
  \matrix (m) [matrix of math nodes,
    nodes in empty cells,nodes={minimum width=5ex,
    minimum height=5ex,outer sep=-5pt},
    column sep=-1ex,row sep=-2ex]{
    	   &&&&&\\
          k     &      &     &     &   & \\
          3     &  \Z &  0 & 2  & \Z & \\
          2     &  0 & 2 &  4  & 2  & \\
          1     &  0 & 2 &  4  & 2  & \\
          0      &  \Z &  0 & 2  & \Z & \\
    \quad\strut &   0  &  1  &  2  & 3  & n \strut \\};
  \draw[thick] (m-1-1.east) -- (m-7-1.east) ;
\draw[thick] (m-7-1.north) -- (m-7-6.north) ;
\end{tikzpicture}
\begin{tikzpicture}
  \matrix (m) [matrix of math nodes,
    nodes in empty cells,nodes={minimum width=5ex,
    minimum height=5ex,outer sep=-5pt},
    column sep=-1ex,row sep=-2ex]{
          &&&&&\\
          k      &      &     &     &   & \\
          3     &  \F &  \F & \F  & \F & \\
          2     &  \F &  \F^2 & \F^2  & \F & \\
          1     &  \F &  \F^2 & \F^2  & \F & \\
          0     &  \F &  \F & \F  & \F & \\
    \quad\strut &   0  &  1  &  2  & 3  & n\strut \\};
  \draw[thick] (m-1-1.east) -- (m-7-1.east) ;
\draw[thick] (m-7-1.north) -- (m-7-6.north) ;
\end{tikzpicture}

\caption{$E_2$-pages of the cohomology spectral sequence $H^n(S^3/\GG;H^k(\T^3;\Z))$ $\Longrightarrow$ $H^{n+k}(\Or;\Z)$ (left) and $H^n(S^3/\GG;H^k(\T^3;\F))$ $\Longrightarrow$ $H^{n+k}(\Or;\F)$ (right). As before, 2 denotes the group $\Z/2$, and $4$ denotes $\Z/4$.}
\label{fig: E2 periodic}
\end{figure}

The computation of $H^*(\Or;\Z)$ now proceeds via these spectral sequences. By $E_n^{*,*}(R )$ we mean the $n^{\rm th}$ page of the spectral sequence $H^*(S^3/\GG;H^*(\T^3;R))\Longrightarrow H^*(\Or;R)$, where $R$ is $\Z$ or $\F$.

\begin{lemma}\label{lem: collapse}
The spectral sequences $H^*(S^3/\GG;H^*(\T^3;R))\Longrightarrow H^*(\Or;R)$ for both $R=\Z$ and $\F$ have no non-trivial differentials.
\end{lemma}

\begin{proof}
The torus $\T^3=\R^3/\Z^3$ has a natural cell structure inherited from the cube, so, with one 0-cell, three 1 and 2-cells and one 3-cell. This cell structure is preserved under the action of $\GG$ and so we may consider the subspaces $Y_n=(X_n\times S^3)/\GG$, $n=0$, 1, 2, where $X_n$ is the $n$-skeleton of $T^3$ with this cell structure. There are then fibrations and inclusions
\[
\begin{array}{ccccc}
\T^3&\lra&\Or&\lra&S^3/\GG\cr
\uparrow&&\uparrow&&||\cr
X_n&\lra&Y_n&\lra&S^3/\GG
\end{array}
\]
induced by the inclusions $i_n\colon X_n\to \T^3$. We write $i_n^*$ both for the induced map in cohomology $H^r(\T^3;R)\to H^r(X_n;R)$, and also for the resulting map of spectral sequences. 

As $i_n^*\colon H^r(\T^3;R)\to H^r(X_n;R)$ is an isomorphism for $r\leqslant n$, and $H^r(X_n;R)=0$ for $r>n$, the $E_2$-page of spectral sequence for $X_n\to Y_n\to S^3/\GG$ with $R$ coefficients is equal to the bottom $n+1$ rows of that for $\T^3\to\Or\to S^3/\GG$ and  $i_n^*$ on the $E_2$-pages is the resulting projection. 

Now suppose there are non-trivial differentials in $E_n^{*,*}(R )$, and suppose $d_m$ is the first, i.e., with smallest $m$. Suppose $x\in E^{a,b}_m(R )$ is an element with $d_m(x)=y\not=0\in E^{a+m, b-m+1}_m(R )$ for minimal $b$. Since $d_m$ is the first non-zero differential, $i^*_{b-1}(y)\not=0$ in $E_m^{*,*}(R)$. However
\[
i^*_{b-1}(y)\ =\ i^*_{b-1}(d_m(x))\ =\ d_m(i^*_{b-1}(x))\ =\ d_m(0)\ =\ 0
\]
since anything in row $b$, that is, in $E_m^{*,b}(R )$, lies in the kernel of $i^*_{b-1}$. This is our contradiction and so there can be no such non-trivial differential $d_m$.
\end{proof}

\begin{proposition}
The cohomology $H^n(\Or;\F)$ of the rotational hull of the cubical tiling of $\R^3$ with $\F$ coefficients is the $\F$ vector space of rank
\[
H^n(\Or;\F)\ =\ \left\{
\begin{array}{ll}
1&n=0\\
2&n=1\\
4&n=2\\
6&n=3\\
4&n=4\\
2&n=5\\
1&n=6\\
0&n>6\,.
\end{array}\right.
\]
\end{proposition}

\begin{proof}
By the previous lemma the spectral sequence $E^{*,*}_*(\F)$ collapses and over the field $\F$ there are no extension problems. The rank of $H^n(\Or;\F)$ can thus be read off the $E_2$-page, counting the ranks of the groups on the diagonal $E_2^{a,n-a}(\F)$ to give the result stated.
\end{proof}

\begin{theorem}
The cohomology $H^n(\Or;\Z)$ of the rotational hull of the cubical tiling of $\R^3$ with integer coefficients is the group

\[H^n(\Or;\Z)\ =\ \left\{
\begin{array}{ll}
\Z&n=0\\
0&n=1\\
\Z/2 \oplus\Z/2              & n=2\\
\Z^2 \oplus \Z/2 \oplus \Z/4 & n=3\\
\Z/2 \oplus \Z/4             & n=4\\
\Z/2 \oplus \Z/2             & n=5\\
\Z                           &n=6\\
0                            &n>6\,.
\end{array}\right.
\]
\end{theorem}

\begin{proof}
The Lemma \ref{lem: collapse} tells us that $E_\infty^{*,*}(\Z)=E_2^{*,*}(\Z)$. Unlike the case of field coefficients there are potential extension problems. However, the case of all extensions being trivial is the only one compatible with the size of the $\F$ coefficient result as stated in the previous proposition: any non-trivial extension would lower the rank of the corresponding $H^n(\Or;\F)$.
\end{proof}


\subsection{Dimension 3: an aperiodic example}\label{sec: Fib3}
We conclude with the computation of the full integer cohomology of the Sturmian decorated cube tiling of Examples \ref{exp: fibcube} and \ref{exp: fibcubecoh}. As in Section \ref{sec: cube}, we use the fibration
\[
\Ot\lra\Or\lra S^3/\GG\,.
\]
Here $\GG$ is still the double cover of the group of symmetries of the cube, but $\Ot$ is shape equivalent to the product of three copies of $W$, the one point union of two circles.  Example \ref{exp: fibcubecoh} sets out the cohomology $H^*(\Ot;\Z)$ and the $\GG$-action.

Computation of group cohomology proceeds analogously to that used in the previous section. As there, we keep track of both the $\Z$ and $\F$ coefficient computations, the latter being used to solve our extension problems at the end. The results for $H^*(B\GG;H^*(\Ot;R))$ are shown in Figure \ref{fig: HBG}.

\begin{figure}
\begin{tikzpicture}
  \matrix (m) [matrix of math nodes,
    nodes in empty cells,nodes={minimum width=5ex,
    minimum height=5ex,outer sep=-5pt},
    column sep=-1ex,row sep=-2ex]{
                &      &     &     &   & \\
          k     &&&&&&&\\
          3     &  \Z^4 &  0 & 2^6  & 0 & 16^2\times48^2 & \cdots &\\
          2     &  0 & 2^3 &  2+4^2  & 2^3  & 0& \cdots&\\
          1     &  0 & 2^2 &  4 ^2 & 2^2  & 0 & \cdots&\\
          0      &  \Z &  0 & 2  & 0 &  48 & \cdots& \\
    \quad\strut &   0  &  1  &  2  & 3  & 4& n \strut \\};
  \draw[thick] (m-1-1.east) -- (m-7-1.east) ;
\draw[thick] (m-7-1.north) -- (m-7-7.north) ;
\end{tikzpicture}
\begin{tikzpicture}
  \matrix (m) [matrix of math nodes,
    nodes in empty cells,nodes={minimum width=5ex,
    minimum height=5ex,outer sep=-5pt},
    column sep=-1ex,row sep=-2ex]{
                &      &     &     &   & &  & \\
          k     &&&&&&&\\
          3     &  \F^4 &  \F^6 & \F^6  & \F^4 & \F^4 & \cdots&\\
          2     &  \F^3 &  \F^6 & \F^6  & \F^3 & \F^3 & \cdots&\\
          1     &  \F^2 &  \F^4 & \F^4  & \F^2 & \F^2 & \cdots&\\
          0     &  \F &  \F & \F  & \F & \F & \cdots&\\
    \quad\strut &   0  &  1  &  2  & 3  &4& n \strut \\};
  \draw[thick] (m-1-1.east) -- (m-7-1.east) ;
\draw[thick] (m-7-1.north) -- (m-7-7.north) ;
\end{tikzpicture}
\caption{Tables of the group cohomologies $H^n(B\GG;H^k(\Ot;\Z))$ (left) and $H^n(B\GG;H^k(\Ot;\F))$ (right). Recall that $H^k(\Ot;R)=0$ for $k>3$. The notation is as before, so $2^2$ denotes the group $\Z/2\times \Z/2$, etc.  All these are 4-periodic, in the sense that  $H^n(B\GG;H^k(\Ot;R))=H^{n+4}(B\GG;H^k(\Ot;R))$ for all $n\geqslant 1$.}
\label{fig: HBG}
\end{figure}


As in the periodic case, we can deduce $H^*(S^3/\GG;M)$ from calculations of $H^*(B\GG;M)$. This allows us to compute the $E_2$-pages of the main spectral sequences in $\Z$ and $\F$ coefficient cohomology, i.e., that for the fibration $\Ot\lra\Or\lra S^3/\GG$. These are as shown in Figure \ref{fig: E2pages}.

\begin{figure}
\begin{tikzpicture}
  \matrix (m) [matrix of math nodes,
    nodes in empty cells,nodes={minimum width=5ex,
    minimum height=5ex,outer sep=-5pt},
    column sep=-1ex,row sep=-2ex]{
                &      &     &     &   & \\
          k     &&&&&&\\
          3     &  \Z^4 &  0 & 2^6  & \Z^4& \\
          2     &  0 & 2^3 &  2+4^2  & 2^3  & \\
          1     &  0 & 2^2 &  4 ^2 & 2^2  & \\
          0      &  \Z &  0 & 2  & \Z &  \\
    \quad\strut &   0  &  1  &  2  & 3  & n \strut \\};
  \draw[thick] (m-1-1.east) -- (m-7-1.east) ;
\draw[thick] (m-7-1.north) -- (m-7-7.north) ;
\end{tikzpicture}
\begin{tikzpicture}
  \matrix (m) [matrix of math nodes,
    nodes in empty cells,nodes={minimum width=5ex,
    minimum height=5ex,outer sep=-5pt},
    column sep=-1ex,row sep=-2ex]{
                &      &     &     &   & \\
          k     &&&&&&\\
          3     &  \F^4 &  \F^6 & \F^6  & \F^4 & \\\
          2     &  \F^3 &  \F^6 & \F^6  & \F^3 & \\
          1     &  \F^2 &  \F^4 & \F^4  & \F^2 & \\
          0     &  \F &  \F & \F  & \F & \\
    \quad\strut &   0  &  1  &  2  & 3  & n \strut \\};
  \draw[thick] (m-1-1.east) -- (m-7-1.east) ;
\draw[thick] (m-7-1.north) -- (m-7-7.north) ;
\end{tikzpicture}
\caption{$E_2$-pages of the spectral sequences $H^n(S^3/\GG;H^k(\Ot;R))\Rightarrow H^{n+k}(\Ot;R)$ for $R=\Z$ and $\F$ respectively. All other rows and columns are zero.}
\label{fig: E2pages}
\end{figure}

\begin{lemma}\label{lem: collapse2}
The spectral sequences $H^*(S^3/\GG;H^*(\Ot;R))\Longrightarrow H^*(\Or;R)$ for both $R=\Z$ and $\F$ have no non-trivial differentials.
\end{lemma}

\begin{proof}
The argument is identical to that used in the proof of  Lemma \ref{lem: collapse}.
\end{proof}

Counting ranks now gives the values of $H^n(\Ot;\F)$, which in turn, as in the previous section, shows there to be no non-trivial extensions in the integer cohomology. We obtain 

\begin{theorem}
The cohomology $H^n(\Or;\Z)$ of the rotational hull of the Sturmian decorated cube tiling with coefficients in $\Z$ and in $\F$ is 
\[
H^n(\Or;\Z)\ =\ \left\{
\begin{array}{ll}
\Z                                   & n=0\\
0                                    & n=1\\
(\Z/2)^3                             & n=2\\
\Z^5 \oplus (\Z/2)^3 \oplus (\Z/4)^2 & n=3\\
(\Z/2)^3 \oplus (\Z/4)^2             & n=4\\
(\Z/2)^9                             & n=5\\
\Z^4                                 & n=6\\
0                                    & n>6
\end{array}\right.
\qquad\qquad\qquad
H^n(\Or;\F)\ =\ \left\{
\begin{array}{ll}
\F      & n=0\\
\F^3    & n=1\\
\F^8    & n=2\\
\F^{15} & n=3\\
\F^{14} & n=4\\
\F^9    & n=5\\
\F^4    & n=6\\
0       & n>6\,.
\end{array}\right.
\]
\end{theorem}


\bibliographystyle{abbrv}
\bibliography{190115HuntonWalton}

\begin{thebibliography}{10}

\bibitem{AaHaOv91}
J.~M. Aarts, C.~L. Hagopian, and L.~G. Oversteegen.
\newblock The orientability of matchbox manifolds.
\newblock {\em Pacific J. Math.}, 150(1):1--12, 1991.

\bibitem{AndPut98}
J.~E. Anderson and I.~F. Putnam.
\newblock Topological invariants for substitution tilings and their associated
  {$C^*$}-algebras.
\newblock {\em Ergodic Theory Dynam. Systems}, 18(3):509--537, 1998.

\bibitem{BaaGri12}
M.~Baake and U.~Grimm.
\newblock On the notions of symmetry and aperiodicity for {D}elone sets.
\newblock {\em Symmetry}, 4(4):566--580, 2012.

\bibitem{BaGr}
M.~Baake and U.~Grimm.
\newblock {\em Aperiodic order. {V}ol. 1}, volume 149 of {\em Encyclopedia of
  Mathematics and its Applications}.
\newblock Cambridge University Press, Cambridge, 2013.
\newblock A mathematical invitation, With a foreword by Roger Penrose.

\bibitem{BSJ91}
M.~Baake, M.~Schlottmann, and P.~D. Jarvis.
\newblock Quasiperiodic tilings with tenfold symmetry and equivalence with
  respect to local derivability.
\newblock {\em J. Phys. A}, 24(19):4637--4654, 1991.

\bibitem{BDHS10}
M.~Barge, B.~Diamond, J.~Hunton, and L.~Sadun.
\newblock Cohomology of substitution tiling spaces.
\newblock {\em Ergodic Theory Dynam. Systems}, 30(6):1607--1627, 2010.

\bibitem{BKS}
M.~Barge, J.~Kellendonk, and S.~Schmieding.
\newblock Maximal equicontinuous factors and cohomology for tiling spaces.
\newblock {\em Fund. Math.}, 218(3):243--268, 2012.

\bibitem{BBG06}
J.~Bellissard, R.~Benedetti, and J.-M. Gambaudo.
\newblock Spaces of tilings, finite telescopic approximations and gap-labeling.
\newblock {\em Comm. Math. Phys.}, 261(1):1--41, 2006.

\bibitem{Bro82}
K.~S. Brown.
\newblock {\em Cohomology of groups}, volume~87 of {\em Graduate Texts in
  Mathematics}.
\newblock Springer-Verlag, New York-Berlin, 1982.

\bibitem{CandelConlon}
A.~Candel and L.~Conlon.
\newblock {\em Foliations. {I}}, volume~23 of {\em Graduate Studies in
  Mathematics}.
\newblock American Mathematical Society, Providence, RI, 2000.

\bibitem{ClHu12}
A.~Clark and J.~Hunton.
\newblock Tiling spaces, codimension one attractors and shape.
\newblock {\em New York J. Math.}, 18:765--796, 2012.

\bibitem{ClHu13}
A.~Clark and S.~Hurder.
\newblock Homogeneous matchbox manifolds.
\newblock {\em Trans. Amer. Math. Soc.}, 365(6):3151--3191, 2013.

\bibitem{ClHuLu14}
A.~Clark, S.~Hurder, and O.~Lukina.
\newblock Shape of matchbox manifolds.
\newblock {\em Indag. Math. (N.S.)}, 25(4):669--712, 2014.

\bibitem{ClHuLu19}
A.~Clark, S.~Hurder, and O.~Lukina.
\newblock Manifold-like matchbox manifolds.
\newblock {\em Proc. Amer. Math. Soc.}, 147(8):3579--3594, 2019.

\bibitem{dWJaJa81}
P.~M. de~Wolff, T.~Janssen, and A.~Janner.
\newblock The superspace groups for incommensurate crystal structures with a
  one-dimensional modulation.
\newblock {\em Acta Cryst. Sect. A}, 37(5):625--636, 1981.

\bibitem{RabFis03}
B.~N. Fisher and D.~A. Rabson.
\newblock Applications of group cohomology to the classification of
  quasicrystal symmetries.
\newblock {\em J. Phys. A}, 36(40):10195--10214, 2003.

\bibitem{FoHuKe02}
A.~Forrest, J.~Hunton, and J.~Kellendonk.
\newblock Topological invariants for projection method patterns.
\newblock {\em Mem. Amer. Math. Soc.}, 159(758):x+120, 2002.

\bibitem{FrWhWh14}
D.~Frettl\"oh, B.~Whitehead, and M.~F. Whittaker.
\newblock Cohomology of the pinwheel tiling.
\newblock {\em J. Aust. Math. Soc.}, 97(2):162--179, 2014.

\bibitem{Gah02}
F.~G\"ahler.
\newblock {Lectures given at workshops {\em Applications of Topology to Physics
  and Biology}, Max-Planck-Institut f\"ur Physik komplexer Systeme, Dresden,
  June 2002, and {\em Aperiodic Order, Dynamical Systems, Operator Algebras and
  Topology}, Victoria, British Columbia, August, 2002.}

\bibitem{GaHuKe13}
F.~G\"ahler, J.~Hunton, and J.~Kellendonk.
\newblock Integral cohomology of rational projection method patterns.
\newblock {\em Algebr. Geom. Topol.}, 13(3):1661--1708, 2013.

\bibitem{GelPro95}
W.~Geller and J.~Propp.
\newblock The projective fundamental group of a {$\bold Z^2$}-shift.
\newblock {\em Ergodic Theory Dynam. Systems}, 15(6):1091--1118, 1995.

\bibitem{Hil86}
H.~Hiller.
\newblock Crystallography and cohomology of groups.
\newblock {\em Amer. Math. Monthly}, 93(10):765--779, 1986.

\bibitem{HunMFO17}
J.~Hunton.
\newblock {\textit{Topological invariants for tilings}, Oberwolfach Reports,
  2017, vol 14(4), 2814-2817}.

\bibitem{Jan88}
T.~Janssen.
\newblock Aperiodic crystals: a contradictio in terminis?
\newblock {\em Phys. Rep.}, 168(2):55--113, 1988.

\bibitem{Ju}
A.~Julien.
\newblock Complexity and cohomology for cut-and-projection tilings.
\newblock {\em Ergodic Theory Dynam. Systems}, 30(2):489--523, 2010.

\bibitem{Kel03}
J.~Kellendonk.
\newblock Pattern-equivariant functions and cohomology.
\newblock {\em J. Phys. A}, 36(21):5765--5772, 2003.

\bibitem{AObook}
J.~Kellendonk, D.~Lenz, and J.~Savinien, editors.
\newblock {\em Mathematics of aperiodic order}, volume 309 of {\em Progress in
  Mathematics}.
\newblock Birkh\"auser/Springer, Basel, 2015.

\bibitem{KelPut06}
J.~Kellendonk and I.~F. Putnam.
\newblock The {R}uelle-{S}ullivan map for actions of {$\Bbb R^n$}.
\newblock {\em Math. Ann.}, 334(3):693--711, 2006.

\bibitem{Mal15}
G.~R. Maloney.
\newblock On substitution tilings of the plane with {$n$}-fold rotational
  symmetry.
\newblock {\em Discrete Math. Theor. Comput. Sci.}, 17(1):395--411, 2015.

\bibitem{MarSeg82}
S.~Marde\v{s}i\'c and J.~Segal.
\newblock {\em Shape theory}, volume~26 of {\em North-Holland Mathematical
  Library}.
\newblock North-Holland Publishing Co., Amsterdam-New York, 1982.
\newblock The inverse system approach.

\bibitem{MarSeg01}
S.~Marde\v{s}i\'{c} and J.~Segal.
\newblock History of shape theory and its application to general topology.
\newblock In {\em Handbook of the history of general topology, {V}ol. 3},
  volume~3 of {\em Hist. Topol.}, pages 1145--1177. Kluwer Acad. Publ.,
  Dordrecht, 2001.

\bibitem{mGMo92}
C.~A. McGibbon and J.~M. M{\o}ller.
\newblock On spaces with the same {$n$}-type for all {$n$}.
\newblock {\em Topology}, 31(1):177--201, 1992.

\bibitem{Mer92}
N.~D. Mermin.
\newblock Erratum: ``{T}he space groups of icosahedral quasicrystals and cubic,
  orthorhombic, monoclinic, and triclinic crystals''.
\newblock {\em Rev. Modern Phys.}, 64(2):635, 1992.

\bibitem{MimTod91}
M.~Mimura and H.~Toda.
\newblock {\em Topology of {L}ie groups. {I}, {II}}, volume~91 of {\em
  Translations of Mathematical Monographs}.
\newblock American Mathematical Society, Providence, RI, 1991.
\newblock Translated from the 1978 Japanese edition by the authors.

\bibitem{MooreSchochet}
C.~C. Moore and C.~L. Schochet.
\newblock {\em Global analysis on foliated spaces}, volume~9 of {\em
  Mathematical Sciences Research Institute Publications}.
\newblock Cambridge University Press, New York, second edition, 2006.

\bibitem{Rad94}
C.~Radin.
\newblock The pinwheel tilings of the plane.
\newblock {\em Ann. of Math. (2)}, 139(3):661--702, 1994.

\bibitem{Ran07}
B.~Rand.
\newblock {\em Pattern-equivariant cohomology of tiling spaces with rotations}.
\newblock ProQuest LLC, Ann Arbor, MI, 2007.
\newblock Thesis (Ph.D.)--The University of Texas at Austin.

\bibitem{Rogers}
J.~W. Rogers, Jr.
\newblock Inducing approximations homotopic to maps between inverse limits.
\newblock {\em Fund. Math.}, 78(3):281--289, 1973.

\bibitem{RoWrMe88}
D.~S. Rokhsar, D.~C. Wright, and N.~D. Mermin.
\newblock Scale equivalence of quasicrystallographic space groups.
\newblock {\em Phys. Rev. B (3)}, 37(14):8145--8149, 1988.

\bibitem{Sad:PEC}
L.~Sadun.
\newblock Pattern-equivariant cohomology with integer coefficients.
\newblock {\em Ergodic Theory Dynam. Systems}, 27(6):1991--1998, 2007.

\bibitem{SadBook08}
L.~Sadun.
\newblock {\em Topology of tiling spaces}, volume~46 of {\em University Lecture
  Series}.
\newblock American Mathematical Society, Providence, RI, 2008.

\bibitem{Sad11}
L.~Sadun.
\newblock Exact regularity and the cohomology of tiling spaces.
\newblock {\em Ergodic Theory Dynam. Systems}, 31(6):1819--1834, 2011.

\bibitem{Ser85}
C.~Series.
\newblock The geometry of {M}arkoff numbers.
\newblock {\em Math. Intelligencer}, 7(3):20--29, 1985.

\bibitem{ShBlGrCa84}
D.~Shechtman, I.~Blech, D.~Gratias, and J.~W. Cahn.
\newblock Metallic phase with long-range orientational order and no
  translational symmetry.
\newblock {\em Phys. Rev. Lett.}, 53:1951--1953, Nov 1984.

\bibitem{Spa66}
E.~H. Spanier.
\newblock {\em Algebraic topology}.
\newblock McGraw-Hill Book Co., New York-Toronto, Ont.-London, 1966.

\bibitem{Sta15}
C.~Starling.
\newblock K-theory of crossed products of tiling {$\rm C^*$}-algebras by
  rotation groups.
\newblock {\em Comm. Math. Phys.}, 334(1):301--311, 2015.

\bibitem{TomZve08}
S.~Tomoda and P.~Zvengrowski.
\newblock Remarks on the cohomology of finite fundamental groups of
  3-manifolds.
\newblock In {\em The {Z}ieschang {G}edenkschrift}, volume~14 of {\em Geom.
  Topol. Monogr.}, pages 519--556. Geom. Topol. Publ., Coventry, 2008.

\bibitem{Wal17pe}
J.~Walton.
\newblock Pattern-equivariant homology.
\newblock {\em Algebr. Geom. Topol.}, 17(3):1323--1373, 2017.

\bibitem{Wal17rot}
J.~J. Walton.
\newblock Cohomology of rotational tiling spaces.
\newblock {\em Bull. Lond. Math. Soc.}, 49(6):1013--1027, 2017.

\end{thebibliography}
\vspace{-0.001cm}
\end{document}